\documentstyle[12pt]{article}  
\def\sq{\hbox {\rlap{$\sqcap$}$\sqcup$}}
\def\sq{\hbox {\rlap{$\sqcap$}$\sqcup$}}
\def\R{ {\rm R \kern -.31cm I \kern .15cm}}
\def\C{ {\rm C \kern -.15cm \vrule width.5pt \kern .12cm}}
\def\Z{ {\rm Z \kern -.27cm \angle \kern .02cm}}
\def\N{ {\rm N \kern -.26cm \vrule width.4pt \kern .10cm}}
\def\1{{\rm 1\mskip-4.5mu l} }
\def\lsim{\raise0.3ex\hbox{$<$\kern-0.75em\raise-1.1ex\hbox{$\sim$}}}
\def\gsim{\raise0.3ex\hbox{$>$\kern-0.75em\raise-1.1ex\hbox{$\sim$}}}
\def\noi{\noindent}

\def\beq{\begin{equation}}   \def\eeq{\end{equation}}
\def\bea{\begin{eqnarray}}  \def\eea{\end{eqnarray}}
\def\nn{\nonumber}
\def\noi{\noindent}
\def\beeq{\begin{eqnarray}} \def\eeeq{\end{eqnarray}}
\newcommand\mysection{\setcounter{equation}{0}\section}

\newcounter{hran}

\begin{document} 
\centerline{\Large\bf Long Range Scattering and Modified } 
 \vskip 3 truemm \centerline{\Large\bf  Wave Operators for the Wave-Schr\"odinger
System III} 

\vskip 0.5 truecm

\centerline{\bf J. Ginibre}
\centerline{Laboratoire de Physique Th\'eorique\footnote{Unit\'e Mixte de
Recherche (CNRS) UMR 8627}}  \centerline{Universit\'e de Paris XI, B\^atiment
210, F-91405 Orsay Cedex, France}
\vskip 3 truemm
\centerline{\bf G. Velo}
\centerline{Dipartimento di Fisica, Universit\`a di Bologna}  \centerline{and INFN, Sezione di
Bologna, Italy}

\vskip 1 truecm

\begin{abstract}
We continue the study of scattering theory for the system consisting of
a Schr\"odinger equation and a wave equation with a Yukawa type
coupling in space dimension 3. In previous papers, we proved the
existence of modified wave operators for that system with no size
restriction on the data and we determined the asymptotic behaviour in
time of solutions in the range of the wave operators, first under a
support condition on the Schr\"odinger asymptotic state and then
without that condition, but for solutions of relatively low regularity.
Here we extend the latter result to the case of more regular solutions. 
\end{abstract}

\vskip 3 truecm
\noi AMS Classification : Primary 35P25. Secondary 35B40, 35Q40, 81U99.  \par
\noi Key words : Long range scattering, modified wave operators, Wave-Schr\"odinger system. 
\vskip 1 truecm

\noindent LPT Orsay 05-12 \par
\noindent February 2005 \par

\newpage
\pagestyle{plain}
\baselineskip 18pt

\mysection{Introduction}
\hspace*{\parindent}
This paper is a sequel to two previous papers with the same title
(\cite{1r} \cite{2r}, hereafter referred to as I and II) where we
studied the theory of scattering and proved the existence of modified
wave operators for the Wave-Schr\"odinger (WS) system in space
dimension 3
$$\hskip 3 truecm \left \{ \begin{array}{l} i\partial_t u = - (1/2) \Delta u - Au \hskip 6.5 truecm (1.1)
\\ \\ \sq A = |u|^2 \hskip 8.9 truecm (1.2)
\end{array}\right .$$

\noi where $u$ and $A$ are respectively a complex valued and a real
valued function defined in space time ${I\hskip-1truemm R}^{3+1}$. We
refer to the introduction of I for general background and references
and we give here only a general overview of the problem. \par

The main result of I and II was the construction of modified wave
operators for the WS system, with no size restriction on the solutions.
That construction basically consists in solving the Cauchy problem for
the WS system with infinite initial time, namely in constructing
solutions with prescribed asymptotic behaviour at infinity in time.
That asymptotic behaviour is imposed in the form of suitable
approximate solutions of the WS system. One then looks for exact
solutions, the difference of which with the given approximate ones
tends to zero at infinity in time in a suitable sense, more precisely
in suitable norms. The approximate solutions are obtained as low order
iterates in an iterative resolution scheme of the WS system. In I and
II we used second order iterates. They are parametrized by data $(u_+,
A_+, \dot{A}_+)$ which play the role of (actually are in simpler cases)
initial data at time zero. Those data constitute the asymptotic state
for the actual solution. \par

An inherent difficulty of the WS system is the difference of
propagation properties of the wave equation and of the Schr\"odinger
equation. Because of that difficulty, we had to impose in I a support
condition on the Fourier transform $Fu_+$ of the Schr\"odinger
asymptotic state $u_+$, saying in effect that $Fu_+$ vanishes in a
neighborhood of the unit sphere, so that $u_+$ generates a solution of
the free Schr\"odinger equation which is asymptotically small in a
neighborhood of the light cone. Such a support condition is unpleasant
because it cannot be satisfied on a dense subspace of any reasonable
space where one hopes to solve the problem, typically with $u$ in
$FH^k$ for $k > 1/2$ ($H^k$ is the standard $L^2$ based Sobolev space).
\par

A progress on that problem was made by Shimomura \cite{13r} \cite{14r}
who was able to construct the wave operators for the WS system without
assuming that support condition, in the framework of a simpler method
proposed earlier by Ozawa \cite{11r} and where the same difficulty
occurs \cite{12r} \cite{15r}. That method however is intrinsically
restricted to the case of small Schr\"odinger data (see \cite{7r} for a
review). The key of that progress consists in using an improved
asymptotic form for the Schr\"odinger function, obtained by adding a
term depending on $(A_+, \dot{A}_+)$ which partly cancels the
contribution of the asymptotic field for $A$ in the Schr\"odinger
equation. \par

Although the method used in I is more complicated than the Ozawa method
(so as to accomodate arbitrarily large data and solutions), it turns
out that the improved asymptotic form of $u$ used in \cite{14r} can be
transposed into the framework of the method of I, thereby allowing to
remove the support condition on $Fu_+$ assumed in I. This was done in
II. However the treatment given in II is restricted to the case of
solutions of relatively low regularity, typically with $u \in FH^k$ with
$1 < k < 3/2$ (the case $k=1$ can also be covered by a simple variant
of the same treatment). The purpose of the present paper is to extend
the same result to the case of more regular solutions, namely with
$k=2$. This is obtained by exploiting the fact that for the
Schr\"odinger equation, one time derivative is homogeneous to two space
derivatives, so that $H^2$ control of the solutions can be achieved by
using only one (time) derivative. That property was extensively used in
similar problems in \cite{4r} \cite{5r} \cite{6r} \cite{15r} (see
\cite{7r} for a review). The methods and results of the present paper
are expected to be relevant for the treatment of the corresponding
problem in the general case for the more complicated
Maxwell-Schr\"odinger system, which was considered in \cite{3r} in the
special case of vanishing asymptotic magnetic field only.\par

In the remaining part of this introduction, we shall briefly review the
method used in I in the modified form used in the present paper. We
refer to Section 2 of I for a more detailed exposition. The main result
of this paper will be stated in semi heuristic terms at the end of this
introduction. The first step in that method consists in eliminating the
wave equation (1.2) by solving it for $A$ and substituting the result
into the Schr\"odinger equation, which then becomes both non linear and
non local in time. One then parametrizes the Schr\"odinger function $u$
in terms of an amplitude $w$ and a phase $\varphi$ and one replaces the
Schr\"odinger equation by an auxiliary system consisting of a transport
equation for the amplitude and a Hamilton-Jacobi equation for the
phase. One solves the Cauchy problem with
prescribed asymptotic behaviour for the auxiliary system, and one
finally reconstructs the solution of the original WS system from that
of the auxiliary system. We now proceed to the technical details. We
restrict our attention to positive time. \par

We first eliminate the wave equation. We define 
$$\omega = (- \Delta)^{1/2} \quad, \quad K(t) = \omega^{-1} \sin \omega t \quad , \quad \dot{K}(t) =
\cos \omega t$$    

\noi and we replace (1.2) by
$$A = A_0 + A_1 (|u|^2) \eqno(1.3)$$

\noi where
$$A_0 = \dot{K}(t) A_+ + K(t) \dot{A}_+ \ , \eqno(1.4)$$
$$A_1 (|u|^2) = - \int_t^{\infty} dt'\ K(t-t') |u(t')|^2 \ . \eqno(1.5)$$

\noi Here $A_0$ is a solution of the free wave equation with initial data $(A_+, \dot{A}_+)$ at time $t
= 0$. The pair $(A_+, \dot{A}_+)$ is the asymptotic state for $A$. \par

We next perform the change of variables mentioned above on $u$. In
connection with the fact that we use time derivatives systematically in
order to perform the necessary estimates, it turns out to be convenient
to perform a different change of variables from that made in I and II.
The main difference lies in the fact that we change $t$ into $1/t$, so
that the auxiliary system will have to be studied for $t$ in a
neighborhood of zero instead of a neighborhood of infinity. The change
of variables from $u$ to $w$ then becomes very similar to the
pseudoconformal inversion (it would reduce to the latter if the phase
$\varphi$ were absent). The unitary group 
$$U(t) = \exp (i(t/2)\Delta ) \eqno(1.6)$$

\noi which solves the free Schr\"odinger equation can be written as
$$U(t) = M(t) \ D(t) \ F \ M(t) \eqno(1.7)$$

\noi where $M(t)$ is the operator of multiplication by the function
$$M(t) = \exp  \left ( i x^2/2t \right ) \ , \eqno(1.8)$$

\noi $F$ is the Fourier transform and $D(t)$ is the dilation operator
$$D(t) = (it)^{-3/2} \ D_0(t) \eqno(1.9)$$

\noi where
$$\left ( D_0(t) f\right )(x) = f(x/t) \ . \eqno(1.10)$$

\noi We parametrize $u$ in terms of an amplitude $w$ and of a real phase $\varphi$ as 
$$u(t) = M(t) \ D(t) \exp [ i \varphi (1/t)] \overline{w}(1/t) \ . \eqno(1.11)$$

\noi Substituting (1.11) into (1.1) yields an evolution equation for $(w , \varphi )$,
namely  
$$\left \{ i \partial_t + (1/2) \Delta - (i/2) (2 \nabla \varphi \cdot \nabla +
\Delta \varphi ) + t^{-1} B + \partial_t \varphi - (1/2) |\nabla \varphi |^2 \right \} w
= 0 \eqno(1.12)$$

\noi where we have expressed $A$ in terms of a new function $B$ by
$$A(t) = t^{-1} \ D_0(t)\ B(1/t) \ . \eqno(1.13)$$

\noi Corresponding to the decomposition (1.3) of $A$, we decompose
$$B = B_0 + B_1 (w , w) \eqno(1.14)$$

\noi where $A_0(t) = t^{-1} D_0(t)B_0(1/t)$ and $A_1(t) =
t^{-1}D_0(t)B_1(1/t)$. One computes easily
$$B_1(w_1,w_2) = \int_1^{\infty} d \nu \ \nu^{-3} \ \omega^{-1} \sin ((\nu - 1) \omega)
D_0(\nu) ({\rm Re} \ \bar{w}_1w_2)(t/\nu ) \ .\eqno(1.15)$$

\indent At this point, we have only one evolution equation (1.12) for
two functions $(w, \varphi )$. We arbitrarily impose a second equation,
namely a Hamilton-Jacobi (or eikonal) equation for the phase $\varphi$,
thereby splitting the equation (1.12) into a system of two equations,
the other one of which being a transport type equation for the
amplitude $w$. For that purpose, we split $B$ into long range and
short range parts as follows. Let $\chi \in {\cal C}^{\infty}
({I\hskip-1truemm R}^{3}, {I\hskip-1truemm R})$, $0 \leq \chi \leq 1$,
$\chi (\xi) = 1$ for $|\xi| \leq 1$, $\chi (\xi) = 0$ for $|\xi | \geq
2$ and let $0 < \beta < 1$. We define
$$B_0 = B_{0L} + B_{0S} \qquad , \qquad B_1 = B_L + B_S \eqno(1.16)$$

\noi where
$$\left \{ \begin{array}{l}  FB_{0L}(t, \xi ) = \chi (\xi t^{\beta})
FB_0(t, \xi ) \ , \\ \\
 FB_L(t, \xi ) = \chi (\xi  t^{\beta})
FB_1(t, \xi ) \ .   \end{array}
\right . \eqno(1.17)$$ 

\noi The parameter $\beta$ will have to satisfy various conditions which will
appear later, all of them compatible with $\beta =  1/3$. \par

We split the equation (1.12) into the following system of two equations.
$$\left \{ \begin{array}{l} i \partial_t w + (1/2)\Delta w =i Q (\nabla \varphi , w)
- t^{-1} (B_{0S} + B_S(w, w)) w \\ \\
 \partial_t \varphi = (1/2) |\nabla \varphi |^2 - t^{-1} \ B_{0L}- t^{-1} \ B_L(w,w)
\end{array} \right . \eqno(1.18)$$ 
\noi where we have defined
$$Q(s, w) = s \cdot \nabla w + (1/2) (\nabla \cdot s) w
\eqno(1.19)$$

\noi for any vector field $s$. The first equation of (1.18) is the
transport type equation for the amplitude $w$, while the second one is
the Hamilton-Jacobi type equation for the phase $\varphi$. Since the
right-hand sides of (1.18) contain $\varphi$ only through its gradient,
we can obtain from (1.18) a closed system for $w$ and $s = \nabla
\varphi$ by taking the gradient of the second equation, namely  
$$\left \{ \begin{array}{l} i\partial_t w + (1/2) \Delta w = i Q (s , w)
- t^{-1} (B_{0S} + B_S(w, w)) w \\ \\
 \partial_t s = s \cdot \nabla s - t^{-1} \nabla  B_{0L} - t^{-1} \nabla  B_L(w,w) \ .\end{array}
\right .  \eqno(1.20)$$ 

 Once the system (1.20) is solved for $(w, s)$, one recovers $\varphi$
easily by integrating the second equation of (1.18) over time. The
system (1.20) will be referred to as the auxiliary system.\par

The construction of the modified wave operators follows the same
pattern as in I and II. The first task is to construct solutions of the
auxiliary system (1.20) with suitably prescribed asymptotic behaviour
at zero, and in particular with $w(t)$ tending to a limit $w_+ =
Fu_+$ as $t \to 0$.  That asymptotic behaviour is imposed in the form
of a suitably chosen pair $(W, \phi )$ and therefore $(W, S)$ with $S =
\nabla \phi$ with $W(t)$ tending to $w_+$ as $t \to 0$. For fixed
$(W,S)$, we make a change of variables in the system (1.18) from $(w,
\varphi )$ to $(q, \psi )$ defined by 
$$(q, \psi) = (w, \varphi ) - (W, \phi )  \eqno(1.21)$$
\noi or equivalently a change of variables in the system (1.20) from $(w, s)$ to $(q, \sigma )$ defined
by   
$$(q, \sigma) = (w, s ) - (W, S) \ , \eqno(1.22)$$

\noi and instead of looking for a solution $(w, s)$ of the system
(1.20) with $(w, s)$ behaving asymptotically as $(W, S)$, we look for a
solution $(q, \sigma )$ of the transformed system with $(q, \sigma )$
(and also $\psi$) tending to zero as $t \to 0$. Performing the change
of variables (1.22) in the auxiliary system (1.20) yields the following
modified auxiliary system for the new variables $(q, \sigma )$
$$\left \{ \begin{array}{l} i\partial_t q +  (1/2) \Delta q = i (Q(s, q) + Q(\sigma,
W)) - t^{-1} B_{0S} \ q  \\
\\ 
-  t^{-1} B_S(w, w) q - t^{-1} \left ( 2 B_S(W,q) + B_S(q,q)
\right ) W - R_1(W,S)\\
\\  
\partial_t \sigma = ( s \cdot \nabla \sigma + \sigma \cdot \nabla S) - t^{-1}\nabla
\left ( 2 B_L(W,q) + B_L(q,q)\right ) - R_2 (W, S) \ ,
 \end{array}
\right . \eqno(1.23)$$

\noi where the remainders $R_1(W, S)$ and $R_2(W,S)$ are defined by 
$$R_1(W,S) = i\partial_t W + (1/2) \Delta W - i Q(S,W) + t^{-1} (B_{0S} + B_S
(W, W)) W \eqno(1.24)$$
$$R_2(W,S) = \partial_tS - S \cdot \nabla S + t^{-1} \nabla B_{0L}+ t^{-1} \nabla B_L(W,W)
\eqno(1.25)$$

\noi and the dependence of the remainders on $B_0$ has been omitted in
the notation. For technical reasons, it is useful to consider also a
partly linearized version of the system (1.23), namely
$$\left \{ \begin{array}{l} i\partial_t q' + (1/2) \Delta q' = i (Q(s, q') + Q(\sigma,
W)) - t^{-1} B_{0S} \ q'  \\
\\ 
- t^{-1} B_S(w, w) q' - t^{-1} \left ( 2 B_S(W,q) + B_S(q,q)
\right ) W - R_1(W,S)\\
\\  
\partial_t \sigma ' = ( s \cdot \nabla \sigma ' + \sigma \cdot \nabla S) - t^{-1}\nabla
\left ( 2 B_L(W,q) + B_L(q,q)\right ) - R_2 (W, S) \ . 
 \end{array}
\right . \eqno(1.26)$$

The construction of solutions $(w,s)$ of the auxiliary system (1.20)
defined for small time and with prescribed behaviour $(W,S)$ at zero
proceeds in two steps. The first step consists in solving the system
(1.23) for $(q, \sigma )$ tending to zero at zero under suitable
boundedness properties of $B_0$ and $(W,S)$ and suitable vanishing
properties of the remainders $R_1(W,S)$ and $R_2(W,S)$ at zero, by a minor
variation of the method used in I. That method consists in first
solving the linearized system (1.26) for $(q' , \sigma ')$ with given
$(q , \sigma )$, and then showing that the map $(q, \sigma ) \to (q',
\sigma ')$ thereby defined has a fixed point, by the use of a
contraction method. The second step consists in constructing $(W, S)$
with $W(t)$ tending to $w_+$ as $t \to 0$ and satisfying the
required boundedness and vanishing properties. This is done by solving the
auxiliary system (1.20) by iteration to second order as in I and then
adding to $W$ an additional term of the same form as that used in
\cite{14r} and in II. The detailed form of $(W, S)$ thereby obtained is
too complicated to be given here and will be given in Section 3 below
(see (3.42)-(3.47)).\par

Once the system (1.20) is solved for $(w, s)$, one can proceed
therefrom to the construction of a solution $(u, A)$ of the original WS
system. One first defines the phases $\varphi$ and $\phi$ such that $s
= \nabla \varphi$ and $S = \nabla \phi$ and one reconstructs $(u, A)$
from $(w, \varphi )$ by (1.11) (1.3) (1.5), thereby obtaining a
solution of the WS system defined for large time and with prescribed
asymptotic behaviour. The modified wave operator for the WS system is
then defined as the map $\Omega : (u_+, A_+, \dot{A}_+) \to (u, A)$.\par

The main result of this paper is the construction of $(u, A)$ from
$(u_+, A_+, \dot{A}_+)$ as described above, together with the
asymptotic properties of $(u, A)$ that follow from that construction.
It will be stated below in full mathematical detail in Proposition 4.1.
We give here a heuristic preview of that result, stripped from most
technicalities. We set $\beta = 1/3$ for definiteness. \\

\noi {\bf Proposition 1.1.} {\it Let $\beta = 1/3$. Let $(u_+, A_+,
\dot{A}_+)$ be such that $w_+ = Fu_+ \in H^{k_+}$ for sufficiently
large $k_+$, that $(A_+, \dot{A}_+)$ be sufficiently regular, and that
$(FA_+, F\dot{A}_+)$ be sufficiently small near $\xi = 0$. Let $(W, S)$
be the approximate solution of the system (1.20) defined by
(3.42)-(3.47). Then \par (1) There exists $\tau = \tau (u_+, A_+,
\dot{A}_+)$, $0 < \tau \leq 1$, such that the auxiliary system (1.20)
has a unique solution ($w, s)$ in a suitable space, defined for $0 < t
\leq \tau$ and such that $(w -W, s - S)$ tends to zero in suitable
norms when $t \to 0$. \par

(2) There exists $\varphi$ and $\phi$ such that $s = \nabla \varphi$,
$S = \nabla \phi$, $\phi (1) = 0$ and such that $\varphi - \phi$
tends to zero in suitable norms when $t \to 0$. Define $(u, A)$ by
(1.11) (1.3) (1.5). Then $(u, A)$ solves the system (1.1) (1.2) for $t
\geq T = \tau^{-1}$ and $(u(t), A(t))$ behaves asymptotically as $(M(t)
D(t) \exp (i \phi (1/t) ) \overline{W}(1/t)$, $A_0 + A_1 (|D(t)
W(1/t)|^2))$ in the sense that the difference tends to zero in suitable
norms (for which each term separately is $O(1)$) when $t \to \infty$.}
\\

The unspecified condition that $(FA_+, F\dot{A}_+)$ be sufficiently
small near $\xi = 0$ can be shown to follow from more intuitive
conditions in $x$-space, consisting of decay conditions at infinity in
space, and, depending on the values of the parameters defining the
relevant function spaces, of some moment conditions on $(A_+,
\dot{A}_+)$.\par

This paper relies on a large amount of material from I and II. In order
to bring out the structure while keeping duplication to a minimum, we
give without proof a shortened logically self-sufficient sequence of
those intermediate results from I and II that are needed, and we provide a
full exposition only for the parts that are new as compared with I and II.
When quoting I or II, we shall use the notation (I.p.q) or (II.p.q) for
equation (p.q) of I or II and Item I.p.q or II.p.q for Item p.q of I or
II, such as Lemma, Proposition, etc.\par

The remaining part of this paper is organized as follows. In Section 2
we collect notation and some estimates of a general nature. In Section
3, we study the Cauchy problem at zero for the auxiliary system
(1.20). We first prove the existence of solutions under suitable
boundedness properties of $(W, S)$ and suitable vanishing properties of
the remainders at zero (Proposition 3.1). We then define $(W,S)$ and
prove that they satisfy the previous properties, which yields the main
result on the Cauchy problem at zero for the auxiliary system (1.20) (Proposition 3.2).
Finally in Section 4, we construct the wave operators for the WS system
(1.1) (1.2) and we derive the asymptotic properties of the solution
$(u, A)$ in their range that follow from the previous results
(Proposition 4.1).

\mysection{Notation and preliminary estimates} 
\hspace*{\parindent}
In this section we introduce some notation and we collect a number of
estimates which will be used throughout this paper. We denote by
$\parallel \cdot \parallel_r$ the norm in $L^r \equiv
L^r({I\hskip-1truemm R}^{3})$ and we define $\delta (r) = 3/2 - 3/r$. 
For any interval $I$ and any Banach space $X$ we denote by ${\cal C}(I,
X)$ the space of strongly continuous functions from $I$ to $X$ and by
$L^{\infty} (I, X)$ the space of measurable essentially bounded
functions from $I$ to $X$. For real numbers $a$ and $b$ we use the
notation $a \vee b = {\rm Max}(a,b)$ and $a\wedge b = {\rm Min} (a,b)$.
In the estimates of solutions of the relevant equations we shall use
the letter C to denote constants, possibly different from an estimate
to the next, depending on various parameters but in general not on the
solutions themselves or on their initial data. We shall use the
notation $C(a_1, a_2, \cdots )$ for estimating functions, also possibly
different from an estimate to the next, depending on suitable norms
$a_1$, $a_2, \cdots$ of the solutions or of their initial data. \par

We shall use the Sobolev spaces $\dot{H}_r^k$ and $H_r^k$ defined for
$- \infty < k < + \infty$, $1 \leq r \leq \infty$ by
$$\dot{H}_r^k = \left \{ u:\parallel u;\dot{H}_r^k\parallel \ \equiv \ \parallel \omega^ku\parallel_r \ <
\infty \right \}$$

\noi and

$$H_r^k = \left \{ u:\parallel u;H_r^k\parallel \ \equiv \ \parallel <\omega>^ku\parallel_r \ <
\infty \right \}$$

\noi where $\omega = (- \Delta)^{1/2}$ and $< \cdot > = (1 + |\cdot |^2)^{1/2}$. The subscript $r$ will
be omitted if $r = 2$ and we shall use
the notation $\parallel w; H^k\parallel\ = |w|_k$.  

We shall look for solutions of the auxiliary system (1.20) such that
$(w, \nabla s) \in {\cal C}(I, H^k \oplus H^{\ell})$ where $I$ is an
interval and where it is understood that $\nabla s \in L^2$ includes
the fact that $s\in L^6$.\par

We shall use extensively the following Sobolev inequalities, stated
here in ${I\hskip-1truemm R}^n$, but to be used only for $n = 3$. \\

\noi {\bf Lemma 2.1.} {\it Let $1 < q$, $r < \infty$, $1 < p \leq
\infty$ and $0 \leq j < k$. If $p = \infty$, assume that $k - j > n/r$.
Let $\sigma$ satisfy $j/k \leq \sigma \leq 1$ and
$$n/p - j = (1 - \sigma )n/q + \sigma (n/r - k) \ .$$

\noi Then the following inequality holds
\beq  
\label{2.1e}
\parallel \omega^j u \parallel_p \ \leq C \parallel u \parallel_q^{1 - \sigma} \ \parallel
\omega^ku \parallel_r^{\sigma} \ .  \eeq} 

\indent The proof follows from the Hardy-Littlewood-Sobolev (HLS)
inequality (\cite{8r}, p.~117) (from the Young inequality if $p =
\infty$), from Paley-Littlewood theory and interpolation.\\

We shall also use extensively the following Leibnitz and commutator
estimates.\\

\noi {\bf Lemma 2.2.} {\it Let $1 < r, r_1, r_3 < \infty$ and
$$1/r = 1/r_1 + 1/r_2 = 1/r_3 + 1/r_4 \ .$$
\noi Then the following estimates hold
\beq  
\label{2.2e}
\parallel \omega^m (uv) \parallel_r \ \leq C \left ( \parallel  \omega^m u \parallel_{r_1} 
 \ \parallel v \parallel_{r_2} + \parallel  \omega^m v \parallel_{r_3} \  
 \parallel u \parallel_{r_4} \right ) \eeq
\noi for $m \geq 0$, and
\beq  
\label{2.3e}
\parallel [\omega^m , u]v\parallel_r \ \leq C \left ( \parallel  \omega^m u \parallel_{r_1} \  
 \parallel v \parallel_{r_2} + \parallel  \omega^{m -1} v \parallel_{r_3} \  
 \parallel \nabla u \parallel_{r_4} \right ) \eeq
\noi for $m \geq 1$, where $[\ ,\ ]$ denotes the commutator.}\\

The proof of those estimates is given in \cite{9r} \cite{10r} with
$\omega$ replaced by $<\omega >$ and follows therefrom by a scaling
argument. \\

We next give some estimates of $B_{0L}$, $B_{0S}$, $B_L$ and $B_S$
defined by (1.16) (1.17). It follows immediately from (1.16) (1.17)
that 
\beq  
\label{2.4e}
\parallel \omega^m B_{0L}\parallel_2 \ \leq \left (2 t^{-\beta}\right )^{m-p}  \parallel  \omega^p
B_{0L} \parallel_2 \  \leq \left (2 t^{- \beta}\right )^{m-p} \parallel \omega^{p} B_0 \parallel_{2} \eeq

\noi for $m \geq p$ and
\beq  
\label{2.5e}
\parallel \omega^m B_{0S}\parallel_2 \ \leq t^{\beta(p-m)}  \parallel  \omega^p B_{0S}
\parallel_2 \  \leq t^{\beta(p-m)} \parallel \omega^{p} B_0 \parallel_{2} \eeq

\noi for $m \leq p$. Similar estimates hold for $B_L$, $B_S$. We shall need estimates of $B_1$ defined by (1.15) and of $\partial_t B_1$. 
From (1.15) it follows  that
\beq
\label{2.6e}
\partial_t B_1 (w_1, w_2) = \widetilde{B}_1 \left ( \partial_t w_1, w_2\right ) + \widetilde{B}_1 \left ( w_1, \partial_t w_2 \right ) 
\eeq

\noi where
\beq
\label{2.7e}
\widetilde{B}_1(w_1,w_2) = \int_1^{\infty} d\nu\ \nu^{-4} \omega^{-1} \sin ((\nu - 1) \omega ) D_0(\nu) \left ( {\rm Re} (\overline{w}_1w_2)\right ) (t/\nu ) \ .
\eeq

\noi It follows from (1.15) (\ref{2.7e}) that
$$\left \{ \begin{array}{l}
\parallel \omega^{m+1} B_{1}(w_1, w_2)\parallel_2 \ \leq I_m \left ( \parallel \omega^{m} (\overline{w}_1
w_2)\parallel_2 \right )\\
\\
\parallel \omega^{m+1} \widetilde{B}_{1}(w_1, w_2)\parallel_2 \ \leq I_{m+1} \left ( \parallel \omega^{m} (\overline{w}_1
w_2)\parallel_2 \right )
\end{array}\right . \eqno(2.8)$$

\noi where $I_m$ is defined by
$$\left ( I_m(f) \right ) (t) = \int_1^{\infty} d\nu \ \nu^{-m-3/2} f(t/\nu ) \ .\eqno(2.9)$$

We finally collect some estimates of the solutions of the free wave equation\break \noindent $\sq A_0 = 0$
with initial data $(A_+, \dot{A}_+)$ at time zero, given by (1.4). \\

\noi {\bf Lemma 2.3.} {\it Let $k \geq 0$. Let $A_+$ and $\dot{A}_+$ satisfy the conditions
$$A_+ , \omega^{-1} \dot{A}_+ \in H^k \qquad , \qquad \nabla^2 A_+, \nabla \dot{A}_+ \in H_1^k \ .\eqno(2.10)$$
Then the following estimate holds~: 
$$\parallel \omega^m A_0 \parallel_r \ \leq \ b_0\ t^{-1+2/r} \qquad {\it for}\  2 \leq r \leq \infty \ ,
\eqno(2.11)$$ 

\noi for $0 \leq m \leq k$ and for all $t > 0$, where $b_0$ depends on $(A_+, \dot{A}_+)$ through the
norms associated with (2.10).} \\

The estimate (2.11) can be expressed in an equivalent form in terms of $B_0$ defined by
(1.13), namely
$$\parallel \omega^m B_0 \parallel_r \ \leq \ b_0\ t^{- m+1/r} \qquad {\it for}\  2 \leq r \leq \infty \ .
\eqno(2.12)$$

\noi We shall also need some estimates on time derivatives of $B_0$. From the fact that the dilation generator $P = t \partial_t + x\cdot \nabla$ satisfies the commutation relation
$$P = D_0(t) \ t\ \partial_t \ D_0 (t)^{-1}\ ,
\eqno(2.13)$$

\noi and from (1.13), it follows that for any non negative integer $j$
$$(P + 1)^j A_0(t) = (-)^j\ t^{-1}\ D_0(t) \ \left ( \left ( t\ \partial_t \right )^j B_0 \right ) (1/t)\ .
\eqno(2.14)$$

\noi On the other hand if $A_0$ is a solution of the wave equation $\sq
A_0 = 0$ with initial data $(A_+, \dot{A}_+)$ at $t=0$, then also
$PA_0$ is a solution of the same equation, with initial data $x\cdot
\nabla A_+$ and $(1 + x \cdot \nabla)\dot{A}_+$. Combining the previous
remarks with Lemma 2.3, we obtain the following corollary.\\

\noi {\bf Corollary 2.1.} {\it Let $A_+$ and $\dot{A}_+$ satisfy the conditions
$$A_+, \omega^{-1} \dot{A}_+ \in L^2 \quad , \quad x \cdot \nabla A_+ ,  \omega^{-1} x\cdot \nabla \dot{A}_+ \in L^2\ .
\eqno(2.15)$$
$$\nabla^2A_+, \nabla \dot{A}_+ \in L^1 \quad , \quad \nabla^2 x \cdot \nabla A_+ , \nabla x\cdot \nabla \dot{A}_+ \in L^1\ .
\eqno(2.16)$$

\noi Then $B_0$ defined by (1.13) satisfies the estimates
$$\parallel \partial_t^j B_0 \parallel_r \leq b_0\ t^{-j+1/r}
\eqno(2.17)$$

\noi for $j = 0, 1$ and for $2 \leq r \leq \infty$, where $b_0$ depends
on $(A_+, \dot{A}_+)$ through the norms associated with (2.15)
(2.16).}\\

Finally from the fact that $\chi \in {\cal C}_0^{\infty}$ and from
obvious scaling properties, it follows that $B_{0S}$ and $B_{0L}$ also
satisfy the estimates (2.12) (2.17), possibly up to an
absolute constant.

\mysection{Cauchy problem at zero for the auxiliary\break\noindent system}
\hspace*{\parindent}
In this section, we solve the Cauchy problem with initial time zero for the auxiliary system (1.20) in the
difference form (1.23). We first solve the system (1.23) for $(q, \sigma )$ tending to zero at zero
under suitable boundedness properties of $(B_0, W, S)$ and suitable vanishing properties of the
remainders $R_1(W,S)$ and $R_2(W,S)$. We then construct $(W,S)$ with $W(t)$ tending to $w_+ = Fu_+$ as
$t \to 0$ and satisfying the required boundedness and vanishing properties. The method closely follows
that of Sections 6 and 7 of I. \par

We first estimate a single solution of the linearized auxiliary system (1.26) at the level of regularity
where we shall eventually solve the auxiliary system (1.20). The following lemma is a variant of Lemma
I.6.1 with $k=2$, where however the second order space derivatives of $q$ are estimated through the use of the time derivative.\\

\noi {\bf Lemma 3.1.} {\it Let $\beta >0$ and $ \ell > 3/2$. Let $\tau
\leq 1$ and $I = (0, \tau]$. Let $B_0 \in {\cal C}^1(I, L^{\infty })$ 
satisfy the estimate (2.17) for $r = \infty$ and $j = 0,1$. Let
$(W,\nabla S) \in {\cal C} (I, H^2 \oplus H^{\ell + 1}) \cap {\cal
C}^1(I, H^1 \oplus H^{\ell - 1})$ and let $(B_0, W, S)$ be such that
$R_1 \in {\cal C}^1(I,L^2)$ and $\nabla R_2 \in {\cal C}(I, H^{\ell})$.
Let $(q, \nabla \sigma )$, $(q', \nabla \sigma ') \in {\cal C}(I, H^2
\oplus H^{\ell})\cap {\cal C}^1 (I, L^2 \oplus H^{\ell-1})$ and let
$(q', \sigma ')$ be a solution of the system (1.26) in I. Assume that
$W$ and $q$ satisfy 

\beq
\label{3.1e}
\mathrel{\mathop {\rm Sup}_{t \in I}} \left ( \parallel W \parallel_{\infty} \ \vee \ |W|_{3/2}\ \vee \ \parallel q
\parallel_{\infty} \ \vee \ |q|_{3/2} \right ) \leq a 
\eeq

\noi for all $t\in I$. Then the following estimates hold for all $t\in I$~:
\beq
\label{3.2e}
\left | \partial_t \parallel q'\parallel_2 \right |\leq \ C \Big \{ a \parallel \nabla
\sigma \parallel_2 \ + \ t^{-1+\beta} \ a^2 \ I_0\left ( \parallel q \parallel_2\right ) \Big \}\  +
\ \parallel R_1 (W, S) \parallel_2  
\eeq
\bea
\label{3.3e}
&&\parallel \Delta q'\parallel_2 \ \leq \ C \Big \{ \parallel \partial_t q'
\parallel_2 \ + \left ( \parallel s \parallel_{\infty}\ + \ \parallel \nabla \cdot s \parallel_3 \right )^2 \parallel q'\parallel_2\ + \ a \parallel \nabla \sigma \parallel_2 \nn \\
&&+ \ t^{-1} \left ( b_0 + a^2\ t^{\beta}\right )  \parallel q'\parallel_2\ + \ t^{-1+\beta} \ a^2 \ I_0\left ( \parallel q \parallel_2\right ) + \ \parallel R_1 (W, S) \parallel_2 \Big \}
 \eea

\noi where $s = S + \sigma$,
\bea
\label{3.4e}
&&\left | \partial_t \parallel \partial_t q'\parallel_2 \right | \leq C \Big \{ \left (\parallel \partial_t s
\parallel_{\infty} \ + \  \parallel \partial_t \nabla \cdot s) \parallel_3\right ) \parallel \nabla q' \parallel_2\ + \ a\parallel \partial_t \nabla \sigma\parallel_2\nn \\
&&+ \left ( \parallel \sigma \parallel_{\infty}\ + \ \parallel \nabla \cdot \sigma \parallel_{3}\right )  \parallel \nabla \partial_t W \parallel_2\ + \ t^{-2}\left ( b_0 + a^2\ t^{\beta}\right )  \parallel q'\parallel_2\nn\\
&&\ + \ t^{-2+\beta} \ a^2 \ I_0\left ( \parallel q \parallel_2\right ) + t^{-1}a \left (\parallel \partial_t W \parallel_3 \ I_0 \left ( \parallel q\parallel_2\right) \ + \ I_0 \left ( \parallel \partial_t W \parallel_3\ \parallel q \parallel_2 \right ) \right )\nn \\
&&\ + \ t^{-1+\beta} \ a^2 \ I_1\left ( \parallel \partial_t q \parallel_2\right ) + t^{-1}a \ I_1\left (\parallel \partial_t W) \parallel_2 \ +\ \parallel \partial_t q\parallel_2\right)  \parallel q'\parallel_3\Big\}\nn \\
&&+ \ \parallel \partial_t R_1 (W, S) \parallel_2 \ ,
\eea
\bea
\label{3.5e}
&&\left | \partial_t \parallel \omega^m \nabla \sigma '\parallel_2 \right | \leq C \Big
\{ \parallel \nabla s\parallel_{\infty} \ \parallel \omega^m \nabla \sigma ' \parallel_2 \
+ \ \parallel \omega^m \nabla s \parallel_2 \ \parallel \nabla \sigma ' \parallel_{\infty}\nn \\
&&+\ \parallel \omega^m \nabla \sigma\parallel_{2} \ \parallel \nabla S
\parallel_{\infty} \ + \ \parallel \sigma \parallel_{\infty} \ \parallel \omega^m \nabla^2
S\parallel_{2} \nn \\
&&+ \ t^{-1- \beta (m+1)} a\ I_0\left ( \parallel q\parallel_2 \right ) \Big \} \ + \ \parallel \omega^m \nabla R_2(W,S)\parallel_2
\eea
\noi for $0 \leq m \leq \ell$, 
$$\left | \partial_t \parallel \nabla \sigma '\parallel_2 \right | \leq C \Big
\{ \parallel \nabla s\parallel_{\infty} \ \parallel \nabla \sigma ' \parallel_2 \
+ \ \parallel \nabla \sigma \parallel_2 \left ( \parallel \nabla S \parallel_{\infty}\ +
\parallel \omega^{3/2} \nabla S \parallel_2 \right ) $$
$$+ \ t^{-1-\beta}\ a\ I_0 \left ( \parallel q \parallel_2\right )  \Big \} \ + \ \parallel \nabla
R_2(W,S)\parallel_2 \ . \eqno(3.5)_0$$
\bea
\label{3.6e}
&& \parallel \partial_t \omega^m \nabla \sigma '\parallel_2  \ \leq C \Big
\{ \parallel s\parallel_{\infty} \ \parallel \omega^m \nabla^2 \sigma ' \parallel_2 \
+ \ \parallel \omega^m \nabla s \parallel_2 \ \parallel \nabla \sigma ' \parallel_{\infty}\nn \\
&&+\ \parallel \omega^m \nabla \sigma\parallel_{2} \ \parallel \nabla S
\parallel_{\infty} \ + \ \parallel \sigma \parallel_{\infty} \ \parallel \omega^m \nabla^2
S\parallel_{2} \nn \\
&&+ \ t^{-1- \beta (m+1)} a\ I_0\left ( \parallel q\parallel_2 \right ) \Big \} \ + \ \parallel \omega^m \nabla R_2(W,S)\parallel_2
\eea

\noi for $0 \leq m \leq \ell - 1$}, 
$$\parallel \partial_t \nabla \sigma '\parallel_2 \ \leq C \Big
\{ \parallel \nabla s\parallel_{\infty} \ \parallel \nabla \sigma ' \parallel_2 \
+ \ \parallel s\parallel_{\infty}\ \parallel \nabla^2 \sigma ' \parallel_2 \ + \ \parallel \nabla \sigma \parallel_2$$
$$\times \left ( \parallel \nabla S \parallel_{\infty}\ +
\parallel \omega^{3/2} \nabla S \parallel_2 \right )  + t^{-1-\beta}\ a\ I_0 \left ( \parallel q \parallel_2\right )  \Big \} \ + \ \parallel \nabla
R_2(W,S)\parallel_2 \ . \eqno(3.6)_0$$
\par \vskip 5 truemm

\noi {\bf Proof}. The estimate (\ref{3.2e}) is essentially (I.6.2)
modified by the change $t \to 1/t$ and simplified by the fact that $q
\in L^{\infty}(I, L^{\infty})$ and is proved in the same way. We next
prove (\ref{3.3e}). From (1.26) and (2.17) we obtain 
$$\parallel \Delta q' \parallel_2 \ \leq 2 \Big\{ \parallel \partial_t q '\parallel_2 \ + \ \parallel Q(s,q')\parallel_2 \ + \ \parallel Q(\sigma, w)\parallel_2 \ + \ C t^{-1}b_0 \parallel q '\parallel_2 $$
\beq
\label{3.7e}
+ \ t^{-1} \parallel B_S(w,w)\parallel_{\infty} \ \parallel q '\parallel_2 \ + \ t^{-1}\parallel B_S(2W+q,q)\parallel_2 \ \parallel W\parallel_{\infty} \ + \parallel R_1\parallel_{2}\Big \} 
\eeq

\noi where $w= W + q$. The terms not containing $q'$ have already been
estimated in the proof of (\ref{3.2e}). We next estimate
\bea
\label{3.8e}
\parallel Q(s, q')\parallel_{2} &\leq& C \left ( \parallel s\parallel_{\infty} \ + \ \parallel \nabla \cdot s\parallel_{3}\right )  \parallel \nabla q'\parallel_2\nn \\
&\leq& C \left ( \parallel s\parallel_{\infty} \ + \ \parallel \nabla \cdot s\parallel_{3} \right ) \parallel q'\parallel_2^{1/2}\parallel \Delta q'\parallel_2^{1/2}\ ,
\eea
\bea
\label{3.9e}
&&\parallel B_S(w,w) \parallel_{\infty} \ \leq \ C\parallel \nabla B_S\parallel_2^{1/2} \ \parallel\nabla^2 B_S \parallel_2^{1/2} \ \leq \ C \ t^{\beta} \parallel \omega^{5/2}B_1(w,w)\parallel_2 \nn \\
&&\leq C \ t^{\beta} \ I_{3/2} \left ( \parallel\omega^{3/2} w \parallel_2 \ \parallel w \parallel_{\infty} \right ) \ \leq \ a^2\ t^{\beta}
\eea

\noi by (\ref{2.5e}) (2.8) and Lemma 2.1. Now (\ref{3.3e})
follows from (\ref{3.2e}) and (\ref{3.7e})-(\ref{3.9e}).\par

We next prove (\ref{3.4e}). Taking the time derivative of the equation
for $q'$ in (1.26), performing a standard $L^2$ norm estimate and using
the fact that the terms in the RHS containing $\partial_t q'$ do not
contribute to that estimate, we obtain
\bea
\label{3.10e}
&&\left | \partial_t \parallel \partial_t q'\parallel_2 \right | \leq  \ \parallel Q(\partial_t s, q')
\parallel_{2} \ + \  \parallel Q( \partial_t \sigma , W) \parallel_2\ + \  \parallel Q (\sigma , \partial_t w)\parallel_2\nn \\
&&+ \ t^{-2} \left \{ \left ( C\ b_0 + \ \parallel B_S(w,w)\parallel_{\infty} \right ) \parallel q' \parallel_{2}\ + \ \parallel B_S(2W + q, q) \parallel_{2}\  \parallel W \parallel_{\infty} \right \}\nn \\
&& + \ t^{-1}\Big \{ \parallel B_S(2W + q,q)\parallel_6\ \parallel \partial_t W \parallel_3 \ + \ \parallel \partial_t B_S(2W + q, q) \parallel_2 \ \parallel W \parallel_{\infty}\nn \\
&& + \ \parallel \partial_t B_S(w,w) \parallel_6 \ \parallel q '\parallel_3\Big \}  \ + \ \parallel \partial_t R_1 \parallel_2\ .
\eea

\noi We next estimate by H\"older and Sobolev inequalities
\beq
\label{3.11e}
\parallel Q(\partial_t s, q')\parallel_2\ \leq\ C \left ( \parallel \partial_t s\parallel_{\infty}\ + \ \parallel \partial_t \nabla \cdot s\parallel_3\right ) \parallel \nabla q'\parallel_2    
\eeq
\bea
\label{3.12e}
\parallel Q(\partial_t \sigma , W)\parallel_2 &\leq& \parallel \partial_t \sigma\parallel_{6}\ \parallel \nabla W \parallel_3\ + \ \parallel \partial_t \nabla \cdot \sigma\parallel_2\  \parallel W\parallel_{\infty}\nn \\
&\leq& C\ a \parallel \partial_t \nabla \sigma \parallel_2    
\eea
\bea
\label{3.13e}
\parallel Q(\sigma , \partial_t  W)\parallel_2 &\leq& \parallel \sigma \parallel_{\infty} \parallel\nabla \partial_t W\parallel_{2}\ + \ \parallel \nabla \cdot \sigma\parallel_3\  \parallel \partial_t W\parallel_{6}\nn \\
&\leq& C\left (  \parallel \sigma \parallel_{\infty}\ + \ \parallel \nabla \cdot \sigma \parallel_3 \right ) \parallel \nabla \partial_t W\parallel_2 \ .   
\eea

\noi Using in addition (\ref{2.5e}) (2.8), we estimate.
\beq
\label{3.14e}
\parallel B_S(2W + q, q)\parallel_2\ \leq\ C\ t^{\beta} \parallel\nabla B_1(2W + q, q)\parallel_2 \ \leq \ C\ a \ t^{\beta}\ I_0 \left ( \parallel q \parallel_2\right )
\eeq
\beq
\label{3.15e}
\parallel B_S(2W + q, q)\parallel_6\ \leq\ C\parallel\nabla B_1(2W + q,q)\parallel_2 \ \leq \ C\ a \ I_0 \left ( \parallel q \parallel_2\right )\ .
\eeq

\noi We next estimate the terms containing $\partial_t B_S$. When
acting on the cut-off factor $\chi$, the time derivative produces a
term 
$$\left ( \partial_t \chi \left ( \xi t^{\beta}\right ) \right ) B_1 = \beta\ t^{-1}\left ( \xi \cdot \nabla \chi \right ) \left ( \xi t^{\beta}\right ) B_1$$

\noi thereby generating an extra factor $t^{-1}$ and a new cut off
field with $\chi$ replaced by $\xi \cdot \nabla \chi$. The
corresponding term is easily estimated in the same way as $B_S$ and
generates terms in the estimates of the same type as those generated by
$\partial_t$ acting on $t^{-1}$. Omitting those terms, we estimate the
remaining contribution of $\partial_t B_S$ as follows (see (\ref{2.6e})-(2.8))
\bea
\label{3.16e}
&&\parallel \partial_t B_S(2W + q, q)\parallel_2\ \leq\ C\Big \{ t^{\beta} \ I_1 \left ( \parallel W + q \parallel_{\infty} \ \parallel \partial_t q \parallel_2 \right )\nn \\
&&+\ I_0 \left ( \parallel \partial_t W \parallel_3\ \parallel q \parallel_2 \right ) \Big \}\nn \\
&&\leq C \left \{ a\ t^{\beta}\ I_1 \left ( \parallel \partial_t q \parallel_2\right ) + I_0\left ( \parallel \partial_t W \parallel_3\ \parallel q \parallel_2\right ) \right \}\ ,
\eea
\bea
\label{3.17e}
&&\parallel \partial_t B_S(w,w )\parallel_6\ \leq\ C\parallel \nabla \partial_t B_S(w,w) \parallel_2 \nn \\
&&\leq C\ I_1  \left ( \parallel w \parallel_{\infty} \ \parallel \partial_t w \parallel_2\right ) \ \leq C\ a \ I_1 \left ( \parallel \partial_t W \parallel_2\ + \ \parallel \partial_t q \parallel_2\right )\ .
\eea

\noi Substituting (\ref{3.9e}) and (\ref{3.11e})-(\ref{3.17e}) into (\ref{3.10e}) yields (\ref{3.4e}). \par

The estimates (\ref{3.5e}) and (3.5)$_0$ are essentially the estimates
(I.6.4) and (I.6.4)$_0$, modified by the change $t \to 1/t$ and
simplified by the fact that $q \in L^{\infty}(I, L^{\infty})$, and are
proved in the same way. The estimates (\ref{3.6e}) and (3.6)$_0$ are
very similar to (\ref{3.5e}) and (3.5)$_0$ and differ therefrom by the
fact that the terms $s\cdot \omega^m\nabla^2\sigma '$ and $s \cdot
\nabla^2 \sigma '$ cannot be integrated by parts when estimating the
$L^2$ norm of the relevant time derivative, whereas they can when
estimating the derivative of the $L^2$ norm.\par\nobreak \hfill $\sq$
\par

We also need to estimate the difference of two solutions of the
linearized auxiliary system (1.26). Those estimates are given by Lemma
I.6.2 in the special case $k = 2$, where they become significantly
simpler. We restate them in the following lemma.\\

\noi {\bf Lemma 3.2.} {\it Let $\beta > 0$ and $\ell > 3/2$. Let $\tau
\leq 1$ and $I = (0, \tau ]$. Let $B_0$, $W$, $S$ satisfy the
assumptions of Lemma 3.1. Let $(q_i , \sigma_i)$ and $(q'_i , \sigma
'_i)$, $i = 1,2$ satisfy the assumptions made on $(q, \sigma )$ and
$(q', \sigma ')$ in Lemma 3.1, and in particular let $(q'_i , \sigma
'_i)$ be solutions of the system (1.26) corresponding to $(q_i ,
\sigma_i)$, $i = 1,2$. Define $(q_{\pm} , \sigma_{\pm}) = (1/2)(q_1 \pm
q_2, \sigma_1 \pm \sigma_2)$ and $(q'_{\pm}, \sigma'_{\pm}) = (1/2)
(q'_1 \pm q'_2, \sigma '_1 \pm \sigma '_2)$. Then the following
estimates hold for all $t\in I$
\beq
\label{3.18e}
\left | \partial_t \parallel  q'_-\parallel_2 \right | \ \leq \ C \left \{ a \parallel \nabla \sigma_- \parallel_2\ +  a^2\ t^{-1+ \beta}\ I_0 \left ( \parallel q_- \parallel_2 \right ) \right \} \ ,
\eeq

$$\left | \partial_t \parallel \nabla \sigma '_-\parallel_2 \right | \ \leq \ C \Big \{ \parallel \nabla s_+\parallel_{\infty}\ \parallel\nabla  \sigma '_- \parallel_2\ +  \ \left ( \parallel \nabla s'_+ \parallel_{\infty} \ + \ \parallel \nabla^2 s'_+ \parallel_3 \right ) \parallel\nabla \sigma_- \parallel_2$$
\beq
\label{3.19e}
+ \ a\ t^{-1- \beta}\ I_0 \left ( \parallel q_- \parallel_2 \right ) \Big \} \ ,
\eeq

\noi where $s_+ = S + \sigma_+$, $s'_+ = S + \sigma '_+$.}\\

We can now solve the Cauchy problem at zero for the auxiliary
system (1.23) under suitable boundedness properties of $(B_0, W, S)$
and suitable vanishing properties of the remainders at zero. This is
the first main result of this section. It corresponds to Proposition
I.6.3, part (2) and to Proposition II.3.1.\\

\noi {\bf Proposition 3.1} {\it Let $0 < \beta < 1$, $\ell > 3/2$ and
$\lambda_0 > 1 \vee \beta (\ell + 1)$. Let $B_0 \in {\cal C}^1((0,1],
L^{\infty})$ satisfy the estimate (2.17) for $r = \infty$ and $j
= 0,1$. Let $(W, \nabla S) \in {\cal C}((0, 1], H^2 \oplus H^{\ell +
1}) \cap {\cal C}^1((0, 1], H^1 \oplus H^{\ell - 1})$, let $R_1 \in
{\cal C}^1((0, 1], L^2)$ and $\nabla R_2 \in {\cal C} ((0, 1],
H^{\ell})$. Let $W$, $S$, $R_1$, $R_2$ satisfy the following estimates
for all $t \in (0,1]$~:
\beq
\label{3.20e}
\parallel W \parallel_{\infty} \ \vee \ |W|_{3/2} \leq a
\eeq
\beq
\label{3.21e}
|\partial_t W|_1 \leq a_1\ t^{-1/2}
\eeq
\beq
\label{3.22e}
\parallel \omega^m \nabla S \parallel_2\ \leq b\left ( |\ell n\ t| + t^{-\beta (m-2)}\right ) \qquad \hbox{\it for $0 \leq m \leq \ell + 1$\ ,}
\eeq
\beq
\label{3.23e}
\parallel \partial_t S \parallel_{\infty} \ \vee \ \parallel \partial_t \nabla S \parallel_3\ \leq b_1\ t^{-1}
\eeq
\beq
\label{3.24e}
\parallel \partial_t^j R_1 \parallel_2\ \leq r_1\ t^{\lambda_0 - 1 - j} \qquad\qquad \hbox{\it for $j=0,1$\ ,}
\eeq
\beq
\label{3.25e}
\parallel \omega^m \nabla R_2 \parallel_2\ \leq r_2\   t^{\lambda_0 - 1 - \beta (m+1)} \qquad \hbox{\it for $0 \leq m \leq \ell$\ .}
\eeq

\noi Then there exists $\tau$, $0 < \tau \leq 1$ and positive constants
$Y_0$, $Y_1$, $Y$, $Z$, depending on $\beta$, $\ell$, $\lambda_0$,
$b_0$, $a$, $a_1$, $b$, $b_1$, $r_1$, $r_2$ such that the auxiliary
system (1.20) has a unique solution $(w, s)$ such that $(w, \nabla s)
\in {\cal C}(I, H^2 \oplus H^{\ell}) \cap {\cal C}^1(I, L^2 \oplus
H^{\ell - 1})$, where $I = (0, \tau ]$, and such that $(q, \sigma )
\equiv (w - W, s-S)$ satisfies the estimates 
\beq
\label{3.26e}
\parallel q \parallel_2\ \leq Y_0\ t^{\lambda_0}
\eeq
\beq
\label{3.27e}
\parallel \partial_t q \parallel_2\ \leq Y_1\ t^{\lambda_0 - 1 }
\eeq
\beq
\label{3.28e}
\parallel \Delta q \parallel_2\ \leq Y\ t^{\lambda_0 - 1}
\eeq
\beq
\label{3.29e}
\parallel \omega^m \nabla \sigma  \parallel_2\ \leq Z \ t^{\lambda_0 - \beta (m+1)} \qquad \hbox{\it for $0 \leq m \leq \ell$\ ,}
\eeq
\beq
\label{3.30e}
\parallel \omega^m \nabla \partial_t \sigma  \parallel_2\ \leq Z \ t^{\lambda_0 - 1 - \beta (m+1)} \qquad \hbox{\it for $0 \leq m \leq \ell - 1$}
\eeq

\noi for all $ t \in I$.}\\

\noi {\bf Sketch of proof.} The proof follows closely that of
Proposition I.6.3, part (2). We first take $\tau$, $0 < \tau \leq 1$,
and $(q, \sigma )$ satisfying the conditions of the proposition and in
particular the estimates (\ref{3.26e})-(\ref{3.30e}) for all $t \in I =
(0, \tau ]$. We then take $t_0$, $0 < t_0 < \tau$, and solve the Cauchy
problem for the linearized system (1.26) with initial condition $(q',
\sigma ') (t_0) = 0$, by the use of Proposition I.6.1 with $k = 2$ and
of Proposition 3.2 in \cite{4r}. Let $(q'_{t_0}, \sigma '_{t_0})$ be the
solution thereby obtained. Using Lemma 3.1, we then show that
$(q'_{t_0}, \sigma '_{t_0})$ satisfies estimates similar to
(\ref{3.26e})-(\ref{3.30e}) with constants $Y'_0$, $Y'_1$, $Y'$, $Z'$,
uniformly in $t_0$ for $t \in [t_0, \tau ]$. Using Lemma 3.2, we take
the limit $t_0 \to 0$ of $(q'_{t_0}, \sigma '_{t_0})$, thereby
obtaining a solution $(q', \sigma ')$ of the system (1.26) in $(0, \tau
]$ satisfying the same estimates in that interval. We finally prove
that for sufficiently small $\tau$, the map $(q, \sigma ) \to (q',
\sigma ')$ thereby defined is a contraction in the norms of Lemma 3.2
on a suitable bounded set defined by the conditions
(\ref{3.26e})-(\ref{3.30e}) for suitably chosen $Y_0$, $Y_1$, $Y$, $Z$.
The abstract arguments are the same as in Proposition I.6.2 and I.6.3,
part (2), and the only difference lies in the estimates of $(q', \sigma
')$ for given $(q, \sigma )$, which now involve a time derivative
instead of space derivatives only. In the remaining parts of this
sketch, we concentrate on the derivation of those estimates. We first
estimate $(q'_{t_0}, \sigma '_{t_0})$ defined above, assuming that $(q,
\sigma )$ satisfies (\ref{3.26e})-(\ref{3.30e}). Omitting the subscript
$t_0$ for brevity, we define
$$\left \{ \begin{array}{l} y'_0 = \ \parallel q'\parallel_2 \ , \qquad y'_1 = \ \parallel \partial_t q'\parallel_2\ , \qquad y' = \ \parallel \Delta q'\parallel_2\\
\\
z'_m = \ \parallel \omega^m \nabla \sigma ' \parallel_2 \qquad \hbox{\rm for $0 \leq m \leq \ell$ \ ,} \end{array} \right . \eqno(3.31)$$

$$\left \{ \begin{array}{l} Y'_0 = \displaystyle{\mathrel{\mathop {\rm Sup}_{t \in I_0 }}}\ t^{- \lambda_0} y'_0 \ , \qquad Y'_1 = \displaystyle{\mathrel{\mathop {\rm Sup}_{t \in I_0 }}}\ t^{1 - \lambda_0} y'_1 , \qquad Y' = \displaystyle{\mathrel{\mathop {\rm Sup}_{t \in I_0 }}}\ t^{1- \lambda_0} y'\\
\\
Z' = \displaystyle{\mathrel{\mathop {\rm Sup}_{0 \leq m \leq \ell }}}\ Z'_m \ , \qquad Z'_m =  \displaystyle{\mathrel{\mathop {\rm Sup}_{t \in I_0 }}}\ t^{\beta (m+1)-\lambda_0} z'_m\ , \end{array} \right . \eqno(3.32)$$

\noi where $I_0 = [t_0, \tau ]$. We first estimate $Y'_0$. From
(\ref{3.2e}) (\ref{3.24e}) (\ref{3.26e}) (\ref{3.29e}) we obtain
$$\left | \partial_t  y'_0\right | \leq C \left ( a \ Z\ t^{\lambda_0 - \beta} + a^2\ Y_0\ t^{\lambda_0 - 1 + \beta}\right ) + r_1\ t^{\lambda_0 - 1} = M(t)\ t^{\lambda_0 - 1} \eqno(3.33)$$

\noi where
$$ M(t) = C \left ( a\ Z\ t^{1-\beta} + a^2\ Y_0 \ t^{\beta}\right ) + r_1 \eqno(3.34)$$

\noi and therefore by integration
$$y'_0 \leq \int_{t_0}^t dt'\ t'^{\lambda_0-1}\ M(t') \leq M(t)\ \lambda_0^{-1} \ t^{\lambda_0} < M(t)\ t^{\lambda_0}$$

\noi so that
$$Y'_0 \leq M (\tau ) \ . \eqno(3.35)$$

\noi We next estimate $Y'$. From (\ref{3.2e}) (\ref{3.3e})
(\ref{3.22e}) (\ref{3.29e}) (3.33) (3.35) we obtain
$$y' \leq C \left \{ y'_1 + M(t) \left ( 1 + b_0 + a^2t^{\beta} + t \left ( b^2|\ell n\ t|^2 + Z^2\ t^{2\lambda_0 - 3 \beta} \right ) \right ) t^{\lambda_0-1}\right \}$$

\noi and therefore
$$Y' \leq C \left \{ Y'_1 + M (\tau ) \left ( 1 + b_0 + a^2 \ \tau^{\beta} + b^2\ \tau |\ell n\ \tau|^2 + Z^2\ \tau^{2\lambda_0 + 1 - 3 \beta}\right ) \right \}\ . \eqno(3.36)$$

\noi We next estimate $Y'_1$. From (\ref{3.4e}) (\ref{3.21e})
(\ref{3.23e}) (\ref{3.24e}) (\ref{3.26e}) (\ref{3.27e}) (\ref{3.30e})
and the definitions, we obtain
$$\left | \partial_t  y'_1\right | \leq C \Big \{ \left ( b_1 \ t^{\lambda_0 - 3/2} + Z\ t^{2\lambda_0 - 3(1 + \beta )/2}\right ) \left ( Y'_0 \ Y'\right )^{1/2}$$
$$+ Z\left ( a\  t^{\lambda_0 - 1 - \beta} + a_1\ t^{\lambda_0 - 1/2 - 3 \beta /2}\right ) + b_0\ Y'_0\ t^{\lambda_0 - 2}$$
$$+ a^2 \left ( Y'_0 + Y_0 + Y_1 \right ) t^{\lambda_0 - 2+ \beta} + a\ a_1\ Y_0\ t^{\lambda_0 - 3/2}$$
$$+ a \left ( a_1\ t^{\lambda_0 - 7/4} + Y_1\ t^{2\lambda_0 - 9/4}\right ) Y'^{3/4}_0 \ Y'^{1/4}\Big \} + r_1\ t^{\lambda_0 - 2}\ .\eqno(3.37)$$

\noi The initial condition for $y'_1$ at $t = t_0$ is estimated by
$$y'_1(t_0) \leq M(t_0)\  t_0^{\lambda_0 - 1} \leq M(t_0)\  t^{\lambda_0-1} \eqno(3.38)$$

\noi for all $t \in I_0$. Integrating (3.37) with initial condition (3.38) and using the definitions yields
$$Y'_1 \leq C \Big \{ \left ( b_1\ \tau^{1/2} + Z\ \tau^{\lambda_0 + (1 - 3 \beta )/2}\right ) \left ( Y'_0\ Y'\right )^{1/2} + a\ Z\ \tau^{1-\beta}$$
$$+ a_1\ Z\ \tau^{3(1-\beta )/2} + b_0\ Y'_0 + a^2\left ( Y'_0 + Y_0 + Y_1 \right ) \tau^{\beta} + a\ a_1 \ Y_0\ \tau^{1/2}$$
$$+ a\left ( a_1 \ \tau^{1/4} + Y_1\ \tau^{\lambda_0 - 1/4}\right ) Y'^{3/4}_0 \ Y'^{1/4} + r_1 \Big \}\ . \eqno(3.39)$$

\noi We next estimate $Z'$. It is sufficient to estimate $Z'_0$ and
$Z'_{\ell}$. The general case follows by interpolation. The estimates
proceed exactly as in I, with minor differences due to the slightly
different assumption made on $S$ and to the simplification produced by
the fact that now $q\in L^{\infty}(I, L^{\infty})$. We obtain (see
(I.6.68) and (I.6.71))
$$Z'_0 \leq C \exp \left ( C \left ( b\ \tau (1 - \ell n\ \tau ) + Z \tau \right ) \right ) \left \{ b\ Z\ \tau (1 - \ell n\ \tau ) + a\ Y_0 + r_2 \right \} \eqno(3.40)$$
$$Z'_{\ell} \leq C \exp \left ( C \left ( b\ \tau (1 - \ell n\ \tau ) + Z \tau \right ) \right ) \left \{ b\ \tau (1 - \ell n\ \tau ) (Z + Z'_0) + \tau\ Z\ Z'_0 + a\ Y_0 + r_2 \right \} \ . \eqno(3.41)$$

\noi From (3.35) (3.36) (3.39) (3.40) (3.41) and from similar easy
estimates of $\parallel \partial_t \omega^m\nabla \sigma '\parallel_2$
that follow from (\ref{3.6e}) (3.6)$_0$ and from the previous estimates
of $\parallel \omega^m \nabla \sigma '\parallel_2$, it follows that
$(q'_{t_0}, \sigma '_{t_0})$ satisfies  estimates similar to
(\ref{3.26e}) -(\ref{3.30e}) with constants $Y'_0$, $Y'_1$, $Y'$ and
$Z'$ satisfying the estimates just derived, for $t \in [t_0 , \tau]$,
uniformly in $t_0$. From there on, the proof proceeds by straightforward
modifications of those of Propositions I.6.2 and I.6.3, part (2). Using
Lemma 3.2, one takes the limit $t_0 \to 0$ of $(q'_{t_0}, \sigma
'_{t_0})$, thereby obtaining a solution $(q', \sigma ')$ of the
linearized system (1.26) in $I = (0, \tau ]$ with the same regularity
and satisfying the same estimates for all $t \in I$. This defines a map
$\phi : (q, \sigma ) \to (q', \sigma ')$ and the previous estimates
prove that this map is bounded in the norm corresponding to
(\ref{3.26e})-(\ref{3.30e}). One finally shows that $\phi$ is a
contraction on the set ${\cal R}$ defined by
(\ref{3.26e})-(\ref{3.30e}) for the norm corresponding to (\ref{3.26e})
and (\ref{3.29e}) with $m= 0$ for suitably chosen $Y_0$, $Y_1$, $Y$,
$Z$ and sufficiently small $\tau$. The proof is a minor variant of that
of Proposition I.6.3, part (2) and will be omitted.\par\nobreak \hfill $\sq$ \par

\noi {\bf Remark 3.1.} Stronger uniqueness results for the system
(1.23) than stated in Proposition 3.1 and stronger uniqueness results
for the system (1.20) regardless on whether $(w, s)$ is an
approximation of some given $(W,S)$ can be obtained by a minor variant
of Proposition I.4.2. For instance if $(w_1, s_1)$ and $(w_2, s_2)$ are
two solutions of the system (1.20) in $I = (0, \tau ]$ satisfying the
regularity properties stated for $(w, s)$ in Proposition 3.1 and the
estimates 
$$\parallel w_i\parallel_{\infty}\ \vee \ |w_i|_{3/2} \leq a$$
$$\parallel \omega^m \nabla s_i\parallel_{2}\ \leq b|\ell n\ t |\qquad \hbox{\rm for $0 \leq m \leq \ell$}$$

\noi for all $t \in I$ or equivalently in a neighborhood of zero and if
$t^{-\beta}\parallel w_1 - w_2\parallel_2$ and $\parallel \nabla (s_1
- s_2)\parallel_2$ tend to zero when $t \to 0$, then $(w_1, s_1) =
(w_2, s_2)$. This can be proved as in Proposition I.4.2, part (3),
starting from a minor variant of Lemma 3.2.\\

We now turn to the construction of approximate solutions $(W,S)$ of the
system (1.20) satisfying the assumptions of Proposition 3.1 and in
particular the estimates (\ref{3.20e})-(\ref{3.25e}) of $W$, $S$ and of
the remainders. As in II, we take for $(W, S)$ the second order
approximate solution of the system (1.20) in an iterative scheme not
taking $B_0$ into account (previously used in I), supplemented by an
additional term in $W$ in order to partly cancel $B_{0S}W$ in $R_1(W,
S)$. Thus we define
$$W = w_0 + w_1 + w_2 \equiv W_1 + w_2\ , \qquad S = s_0 + s_1 \eqno(3.42)$$

\noi where (up to the change of $t$ into $1/t$) $w_0$, $w_1$, $s_0$, $s_1$ are the same as in I except for a simplification of $w_1$ namely
$$w_0 = U(t) w_+ \eqno(3.43)$$
$$s_0(t) = \int_t^1 dt'\ t'^{-1} \nabla B_L \left ( w_0(t') , w_0 (t')\right ) \eqno(3.44)$$
$$w_1(t) = \int_0^t dt'\ Q \left ( s_0(t'), w_0 (t')\right ) \eqno(3.45)$$
$$s_1(t) = \int_0^t dt'\left ( s_0(t') \cdot \nabla s_0 (t') - 2 t'^{-1} \nabla B_L \left ( w_0 (t'), w_1 (t')\right )\right ) \eqno(3.46)$$

\noi while $w_2$ is the same as in II, namely
$$w_2 = h\ w_0 \quad , \qquad h = - 2 t^{-1}\ \Delta^{-1}\ B_{0S}\ . \eqno(3.47)$$

\noi With that choice, the remainders become
$$R_1(W, S) = (1/2) \Delta w_1 + i \left ( \partial_t h\right ) w_0 + (\nabla h)\cdot \nabla w_0 - i Q \left ( s_0, w_1 + w_2 \right )$$
$$- i Q\left ( s_1, W\right ) + t^{-1}\ B_{0S} \left ( w_1 + w_2 \right ) + t^{-1}\ B_S (W, W) W \eqno(3.48)$$
$$R_2(W, S) = R_{20} (W,S) + R_{2\nu } (W, S) \eqno(3.49)$$

\noi where
$$R_{20}(W,S) = - \left ( s_0 \cdot \nabla s_1 + s_1 \cdot \nabla s_0 + s_1 \cdot \nabla s_1 \right ) + t^{-1} \nabla B_L \left ( w_1, w_1\right ) \eqno(3.50)$$
$$R_{2\nu}(W, S) = t^{-1}\nabla B_{0L} + t^{-1} \nabla B_L \left ( W + W_1 , w_2 \right ) \ . \eqno(3.51)$$ 

\noi The terms not involving $h$ or $w_2$ have already been estimated
in I. Up to appropriate minor changes and additions, the following
lemma is basically Lemma I.7.1 or Lemma II.3.2.\\

\noi {\bf Lemma 3.3.} {\it Let $0 < \beta < 1$, $k_+ \geq 3$, $w_+ \in H^{k_+}$ and $a_+ =
|w_+|_{k_+}$. Then the following estimates hold for all $t$, $0 < t \leq 1$~:}
$$|\partial_t w_0|_{k_+-2} \ \leq |w_0|_{k_+} = a_+ \eqno(3.52)$$ 
$$\parallel \omega^m \ s_0\parallel_2 \ \leq C\ a_+^2\left ( |\ell n \ t| + t^{-\beta (m-k_+)}\right ) \qquad \hbox{\it for} \ m \geq 0\ , \eqno(3.53)$$
$$\parallel \omega^m \ \partial_t  s_0\parallel_2 \ \leq C\ a_+^2\ t^{-1-\beta (m-k_+)_+} \qquad \hbox{\it for} \ m \geq 0\ , \eqno(3.54)$$
$$|w_1|_{k_+ - 1} \leq C \ a_+^3\ t(1 - \ell n\ t) \eqno(3.55)$$ 
$$|\partial_t \ w_1 |_{k_+-1} \leq C\ a_+^3|\ell n \ t| \eqno(3.56)$$
$$\parallel \omega^m \ s_1\parallel_2 \ \leq C\ a_+^4 \ t\left \{  (1
- \ell n \ t)^2  + (1 - \ell n \ t) t^{-\beta (m +1- k_+)}\right \} \eqno(3.57)$$
\hskip 8 truecm  $\hbox{\it for}\ m \geq 0 , \ \beta (m+1 - k_+) < 1 \ ,$ 
$$\parallel \omega^m \ \partial_t  s_1\parallel_2 \ \leq C\ a_+^4 \left \{  (1
- \ell n \ t)^2  + (1 - \ell n \ t) t^{-\beta (m +1- k_+)}\right \} \eqno(3.58)$$
\hskip 8 truecm  $\hbox{\it for}\ m \geq 0 , \ \beta (m+1 - k_+) < 1 \ ,$ 
$$\parallel \omega^m \ R_{20}(W,S)\parallel_2\   \leq C(a_+) t
\left \{ (1 - \ell n \ t)^3 + (1 - \ell n\ t)^2  t^{-\beta (m +2- k_+)}\right \}\eqno(3.59)$$
\hskip 8 truecm  $ \hbox{\it for}\ m \geq 0\ , \beta (m+2 - k_+) < 1 \ . $\par \vskip 5 truemm

In order to complete the estimates of $W$ and of the remainders, we
need some estimates of $B_{0L}$ and of $h$. Those estimates are the
same as in II, supplemented by an estimate of $\partial_t^2h$. They
require some restrictions on the behaviour of $(FA_+, F \dot{A}_+)$ at
$\xi = 0$, which we impose as in II in a dilation homogeneous way in
terms of quantities which have the same scaling properties as
$\parallel A_+; \dot{H}^{-3/2-\mu}\parallel$ and $\parallel \dot{A}_+;
\dot{H}^{-5/2-\mu}\parallel$ for some $\mu \in (-1, 1)$. The following
lemma is a minor extension of Lemma II.3.4 and its proof will be
omitted.\\

\noi {\bf Lemma 3.4.} {\it Let $0 < \beta < 1$ and $-1 < \mu < 1$. Let $(A_+, \dot{A}_+)$ satisfy the
conditions
$$A_+, xA_+, x^2 A_+ \in H^{1} \qquad , \qquad \dot{A}_+, x\dot{A}_+, x^2 \dot{A}_+  \in L^2 \eqno(3.60)$$
$$ x^2 A_+ \in \dot{H}^{1/2-\mu} \ , \ xA_+ ,  x^2 \dot{A}_+ \in \dot{H}^{-1/2 - \mu} \ , \
A_+, x\dot{A}_+ \in \dot{H}^{-3/2 - \mu} \ , \dot{A}_+ \in \dot{H}^{-5/2 - \mu} \ . \eqno(3.61)_{\mu}$$

\noi Let $B_{0L}$ and $h$ be defined by (1.17) and (3.47). Then the following estimates hold~:
$$\parallel \omega^m B_{0L} \parallel_2 \ \leq C\ t^{2 + \mu - \beta (m+3/2 + \mu )} \left (
\parallel A_+ ; \dot{H}^{-3/2-\mu} \parallel \  + \ \parallel \dot{A}_+ ; \dot{H}^{-5/2 - \mu} \parallel
\right ) \eqno(3.62)$$

\noi for all $m \geq 0$,
$$\parallel \omega^m \ \partial_t^j h \parallel_2 \ \leq C\ t^{3/2 -m-j} 1_m \sum_{0 \leq j' \leq j} \left ( \parallel |x|^{j'} A_+; \dot{H}^{\rho + j'} \parallel \ + \ \parallel |x|^{j'} \dot{A}_+; \dot{H}^{\rho + j' -
1} \parallel \right )\eqno(3.63)$$

\noi for $j = 0, 1, 2$ and $m \leq 3 - j$, where
$$1_m =  1 \vee t^{(1- \beta )(m-1/2 + \mu)} \eqno(3.64)$$

\noi and 
$$\rho = (m-2) \vee (-3/2 - \mu) \ , \eqno(3.65)$$
$$\parallel h \parallel_{\infty} \ \leq C(A_+, \dot{A}_+) \eqno(3.66)$$

\noi where the constant depends on $(A_+, \dot{A}_+)$ through the norms in (3.60)
(3.61)$_{\mu}$ not involving $x$.} \\

We can now state the final result on the Cauchy problem at zero for the auxiliary system (1.20).\\

\noi {\bf Proposition 3.2.} {\it Let $0 < \beta < 3/5$, $\ell > 3/2$
and $1 \vee \beta (\ell + 1) < \lambda_0 < 3/2$. Let $k_+$ and $\mu$
satisfy
$$k_+ \geq \lambda_0 + 2 \qquad , \qquad \beta \left ( k_+ + 1 \right ) \geq \lambda_0 \ , \eqno(3.67)$$
$$\mu \geq (1 - \beta )^{-1} \left \{ \lambda_0 - 1 - \beta /2 + (3 \beta - 1) \vee 0 \right \} \ (> - 1/4)\ . \eqno(3.68)$$

\noi Let $w_+ \in H^{k_+}$, and let $(A_+, \dot{A}_+)$ satisfy
(2.16) (3.60) (3.61)$_{\mu}$. Let $(W, S)$ be defined by
(3.42)-(3.47). Then there exists $\tau$, $0 < \tau \leq 1$ such that
the auxiliary system (1.20) has a unique solution $(w, s)$ such that
$(w, \nabla s) \in {\cal C} (I, H^2 \oplus H^{\ell}) \cap {\cal C}^1(I,
L^2 \oplus H^{\ell - 1})$, where $I = (0 , \tau ]$, and satisfying
$$\parallel w(t) - W(t) \parallel_2\ \leq C \left ( a_+, A_+, \dot{A}_+\right ) t^{\lambda_0} \eqno(3.69)$$ 
$$\parallel \partial_t (w(t) - W(t)) \parallel_2\ \vee \ \parallel \Delta (w(t) - W(t))\parallel_2\ \leq C \left ( a_+, A_+, \dot{A}_+\right ) t^{\lambda_0-1} \eqno(3.70)$$ 
$$\parallel \omega^m (s(t) - S(t)) \parallel_2\ \leq C \left ( a_+, A_+, \dot{A}_+\right ) t^{\lambda_0-\beta m} \ \hbox{\it for $0 \leq m \leq \ell + 1$} \eqno(3.71)$$ 
$$\parallel \omega^m \ \partial_t (s(t) - S(t)) \parallel_2\ \leq C \left ( a_+, A_+, \dot{A}_+\right ) t^{\lambda_0-1 - \beta m} \ \hbox{\it for $0 \leq m \leq \ell$} \eqno(3.72)$$

\noi for all $t \in (0, \tau ]$, where $a_+ = |w_+|_{k_+}$, and the
constants $C(a_+, A_+, \dot{A}_+)$ depend on $(A_+, \dot{A}_+)$ through
the norms associated with (2.16) (3.60) (3.61)$_{\mu}$.}\\

\noi {\bf Proof.} Proposition 3.2 follows from Proposition 3.1 and from
the fact that $(W, S)$ defined by (3.42)-(3.47) satisfies the
assumptions of the latter under the assumptions made here, and in
particular satisfies the estimates (\ref{3.20e})-(\ref{3.25e}). The
regularity properties are easily seen to hold and we concentrate on the
estimates. The estimates (\ref{3.20e}) (\ref{3.21e}) as regards $W_1 =
w_0 + w_1$ and the estimates (\ref{3.22e}) (\ref{3.23e}) follow
immediately from Lemma 3.3. We estimate $w_2$ by Lemma 3.4 and possibly
Lemmas 2.1 and 2.2, namely 
$$\parallel w_2\parallel_2\ \leq \ \parallel h \parallel_2\ \parallel w_0 \parallel_{\infty}\ \leq C\ t^{3/2} \ 1_0\ , \eqno(3.73)$$ 
$$\parallel \omega^m\ w_2\parallel_2\ \leq \ C\left ( \parallel \omega^m h \parallel_2\ \parallel w_0 \parallel_{\infty}\ + \ \parallel h \parallel_r\ \parallel \omega^m\ w_0 \parallel_{3/\delta}\right )$$

\noi by Lemma 2.2 with $2 < r < \infty$, so that $0 < \delta \equiv 3/2 - 3/r < 3/2$, 
$$\cdots \leq C \ \left ( \parallel \omega^m h \parallel_2\ \parallel w_0 \parallel_{\infty}\ + \ \parallel \omega^{\delta} h\parallel_2\ \parallel \omega^{m + 3/2 - \delta} w_0 \parallel_2\right )$$

\noi by Lemma 2.1, 
$$\cdots \leq C \ \left ( t^{3/2-m} \ 1_m + t^{3/2- \delta }\ 1_{\delta}\right ) \leq C\ t^{3/2-m} \ 1_m \eqno(3.74)$$

\noi for $0 < m < 3$ with $k_+ \geq 3$ by Lemma 3.4, provided $m - 3/2 \leq \delta \leq m$. One can take for instance $\delta = m/2$.\par

We next estimate
$$\parallel  \partial_t\ w_2 \parallel _2\ \leq \ \parallel  \partial_t  h \parallel _2 \ \parallel w_0\parallel _{\infty}\ + \ \parallel  h \parallel_3 \ \parallel  \Delta w_0\parallel_6 \ \leq\ C\ t^{1/2}\ 1_0\ , \eqno(3.75)$$
$$\parallel  \partial_t\ \nabla w_2 \parallel _2\ \leq \ \parallel  \partial_t \nabla h \parallel _2 \ \parallel w_0\parallel _{\infty}\ + \ \parallel \partial_t  h \parallel_2 \ \parallel  \nabla w_0\parallel_{\infty}$$
$$+\  \parallel  \nabla h \parallel_3 \ \parallel  \Delta w_0\parallel_6 \ + \ \parallel h\parallel _{\infty}\ \parallel  \nabla \Delta w_0\parallel_2\ \leq\ C\ t^{-1/2} \eqno(3.76)$$

\noi which completes the proof of (\ref{3.20e}) (\ref{3.21e}). We now
turn to the remainders and to the proof of (\ref{3.24e}) (\ref{3.25e}).
The $L^2$ estimate of $R_1$ has already been proved in II but we sketch
its proof again because (i) it is simpler in the present case and (ii)
it is the starting point to prove the estimate of $\partial_tR_1$,
which is new. We estimate
$$\parallel R_1\parallel_2\ \leq\ \parallel \Delta w_1 \parallel_2\ + \ \parallel \partial_t h \parallel_2 \ \parallel w_0\parallel_{\infty}\ + \ \parallel \nabla h \parallel_2\ \parallel \nabla w_0\parallel_{\infty}$$
$$+ \left ( \parallel s_0\parallel_{\infty} \ + \ C \parallel \nabla \cdot s_0 \parallel_3 \right ) \left ( \parallel \nabla w_1 \parallel_2 \ + \ \parallel \nabla w_2 \parallel_2 \right )$$
$$+ \left ( \parallel s_1\parallel_{\infty} \ + \ C \parallel \nabla \cdot s_1 \parallel_3 \right ) \parallel \nabla W\parallel_{2} \ + \ t^{-1} \parallel B_0\parallel_{2} \ \parallel w_1 \parallel_{\infty}$$
$$+ \ t^{-1} \parallel B_0\parallel_{\infty} \parallel w_2 \parallel_2 \ + \ t^{-1} \parallel W\parallel_{\infty}\ \parallel B_S(W,W)\parallel_{2}\ . \eqno(3.77)$$

\noi We shall need a slightly more general estimate of $B_S(W,W)$ than required for (3.77). We separate
$$B_S(W,W) = B_S(w_0, w_0) + B_S(w_1, w_0 + W_1) + B_S(w_2, W + W_1) \eqno(3.78)$$

\noi and we estimate
$$\parallel \omega^m \ B_S(w_0, w_0)\parallel\ \leq C\ t^{\beta (k_+ +1-m)}\ I_{k_+} \left ( \parallel \omega^{k_+} w_0 \parallel_2\ \parallel w_0 \parallel_{\infty} \right )$$
$$\leq  C \ a_+^2\ t^{\beta (k_+ + 1 - m)} \eqno(3.79)$$

\noi for $0 \leq m \leq k_+ + 1$, by (\ref{2.5e}) (2.8) and Lemma 2.2,
$$\parallel \omega^m \ B_S(w_1, w_0 + W_1)\parallel_2\ \leq C\ t^{\beta (k_+-m)}\ I_{k_+-1} \Big ( \parallel \omega^{k_+-1} w_1 \parallel_2$$
$$\times \left ( \parallel w_0 \parallel_{\infty} \ + \ \parallel w_1 \parallel_{\infty}\right ) + \parallel \omega^{k_+} w_0 \parallel_{2} \ \parallel w_1 \parallel_{3}\Big )$$
$$ \leq  C (a_+)t^{\beta (k_+  - m)+1} (1 - \ell n\ t)\eqno(3.80)$$

\noi for $0 \leq m \leq k_+$, by (\ref{2.5e}) (2.8), Lemma 2.2 and (3.55),
$$\parallel \omega^m \ B_S(w_2, W + W_1)\parallel_2\ \leq C\ t^{\beta (1-m)}\ I_0 \left ( \parallel w_2 \parallel_2 \left ( \parallel w_2 \parallel_{\infty} \ + \ \parallel W_1 \parallel_{\infty}\right ) \right )$$
$$ \leq  C \ t^{\beta (1 - m)+3/2} 1_0\eqno(3.81)$$

\noi for $0 \leq m \leq 1$, by (\ref{2.5e}) (2.8) (3.73), so that
$$\parallel \omega^m \ B_S(W, W)\parallel_2\ \leq C\left ( t^{\beta (k_+ + 1-m)} + t^{\beta (1- m) + 3/2} 1_0  \right ) \eqno(3.82)$$

\noi for $0 \leq m \leq 1$. From (3.77) (2.17), from Lemmas 3.3 and 3.4 and from (3.82) with $m = 0$, it follows that
$$\parallel R_1 \parallel_2\ \leq C \left \{ t^{1/2}\left ( |\ell n\ t | + 1_0\right ) + t^{\beta (k_+ + 1)-1}\right \} \eqno(3.83)$$

\noi so that $R_1$ satisfies (3.24) with $j = 0$ provided $\lambda_0 < 3/2$, $\lambda_0 \leq \beta (k_+ + 1)$ and
$$\lambda_0 \leq 1 + \mu - \beta (\mu - 1/2) \eqno(3.84)$$

\noi which is part of (3.68). We next estimate $\partial_t R_1$ in
$L^2$. The estimate of $\partial_tR_1$ is obtained from (3.77) by
inserting an additional $\partial_t$ at all possible places and
possibly changing the exponents in the application of the H\"older
inequality. If no such change is required, it follows from
(2.17) and from Lemmas 3.3 and 3.4 that the estimate gets worse
by at most one power of $t$ so that all such terms satisfy the $j = 1$
case of (\ref{3.24e}). It is therefore sufficient to consider only the
terms where the H\"older exponents have to be changed. This comes about
mostly as a consequence of the fact that $\partial_t w_0 \sim \Delta
w_0$ has limited regularity, and in particular is in general not in
$L^{\infty}$ under the assumption made on $k_+$. Thus we estimate
$$\parallel \partial_t R_1\parallel_2\ \leq \ \parallel \partial_t h\parallel_3\ \parallel \partial_t  w_0\parallel_6\ + \ C \parallel \nabla h \parallel_r \ \left | \nabla \partial_t w_0\right |_{k_+-3}$$
$$+ \ t^{-1} \left \{ \parallel \partial_t W\parallel_6 \ \parallel B_S(W,W)\parallel_3\ + \ 2\parallel W\parallel_{\infty}\ \parallel \widetilde{B}_S(\partial_t w_0, W\parallel_2 \right \}$$
$$+ \ \hbox{\rm other terms}\ , \eqno(3.85)$$

\noi where $3/r = k_+ - 3$ for $3 < k_+ < 9/2$, the only dangerous
case, and $r=2$ for $k_+ > 9/2$, and where $\widetilde{B}$ is defined by (\ref{2.7e}). We next estimate
$$\parallel \partial_t h \parallel_3 \ \parallel \partial_t w_0\parallel_6 \ \leq\ C \parallel \omega^{1/2}\ \partial_t h\parallel_2\ \parallel \nabla \Delta w_0\parallel_2$$
$$ \leq C\ 1_{1/2} \leq C\ t^{\lambda_0 -2} \eqno(3.86)$$

\noi for $k_+ \geq 3$, 
$$\parallel \nabla h \parallel_r \ \left | \nabla \partial_t w_0 \right |_{k_+-3} \ \leq\ \parallel \omega^{11/2-k_+}\  h\parallel_2\ |w_0|_{k_+}$$
$$ \leq C\ t^{k_+ - 4} \leq C\ t^{\lambda_0 -2} \eqno(3.87)$$

\noi for $\lambda_0 + 2 \leq k_+ < 9/2$,
$$t^{-1} \parallel \partial_t W\parallel_6\ \parallel B_S(W,W)\parallel_3\ \leq\ C\ t^{-1} \parallel \nabla \partial_t W\parallel_2\ \parallel \omega^{1/2} B_S(W,W)\parallel_2$$
$$\leq C\ t^{-3/2} \left ( t^{\beta (k_+ + 1/2)} + t^{(\beta + 3)/2}\ 1_0 \right ) \leq C\ t^{\lambda_0-2} \eqno(3.88)$$

\noi for $\lambda_0 \geq \beta (k_+ + 1)$, by (3.76) and (3.82) with $m = 1/2$,  
$$t^{-1} \parallel \widetilde{B}_S(\partial_t  w_0, w_0)\parallel_2\ \leq C\ t^{-1+\beta (k_+ - 1)} \ I_{k_+-1} \Big ( \parallel \omega^{k_+} w_0\parallel_2\ \parallel w_0\parallel_{\infty}$$
$$+\ \parallel \Delta w_0\parallel_2\ \parallel \omega^{k_+-2} w_0\parallel_{\infty}\Big ) \ \leq C\ t^{-1+\beta (k_+-1)}\leq C\ t^{\lambda_0-2} \eqno(3.89)$$

\noi for $k_+ \geq 2$ and $\lambda_0 \leq \beta (k_+ - 1) + 1$. The
last condition for $k_+ \geq 2$ follows from $\lambda_0 \leq 3/2$ if
$\beta \geq 1/2$ and from $\lambda_0 \leq \beta (k_+ + 1)$ if $\beta
\leq 1/2$. 
$$t^{-1} \parallel \widetilde{B}_S(\partial_t  w_0, w_1 + w_2)\parallel_2\ \leq C\ t^{-1+\beta /2} \ I_{1/2} \left ( \parallel \partial_t w_0\parallel_6\ \parallel w_1 + w_2\parallel_{2}\right ) $$
$$\leq C\ t^{\beta /2}\left ( 1 - \ell n\ t + t^{1/2}\ 1_0\right ) \leq C\ t^{\lambda_0-2} \eqno(3.90)$$

\noi where the last inequality is largely satisfied for $\lambda_0 \leq
3/2$. From (3.85)-(3.90) and the fact that the other terms in (3.85)
are correctly estimated by the previous remark, it follows that $R_1$
satisfies the estimate (\ref{3.24e}) for $j=1$.\par 

The estimate (\ref{3.25e}) of $R_2$ is the same as in II and we briefly recall its
proof for completeness. From (3.59) it follows that $R_{20}$ satisfies
(\ref{3.25e}) provided $\lambda_0 < 2 + \beta$. We next consider
$R_{2\nu}$ defined by (3.51). For $m \geq 0$, we estimate
$$t^{-1} \parallel \omega^{m+1} B_{0L}\parallel_2\ \leq \ C\ t^{1 + \mu - \beta (m + 5/2 + \mu )}\leq C\ t^{\lambda_0 - 1 - \beta m} \eqno(3.91)$$

\noi by (3.62) provided
$$\beta (5/2 + \mu ) \leq 2 + \mu - \lambda_0 \eqno(3.92)$$

\noi which is half of (3.68), and 
$$t^{-1} \parallel \omega^{m+1} B_{L}(W+ W_1, w_2)\parallel_2
\leq t^{-1- \beta m}\ I_0 \left ( \parallel W + W_1\parallel_{\infty}\ \parallel w_2 \parallel_2 \right )$$
$$ \leq C \ t^{1/2 - \beta m}\ 1_0
\leq C\ t^{\lambda_0 - 1 - \beta m} \eqno(3.93)$$

\noi under the condition (3.84) which is the other half of (3.68). This
completes the proof of the estimates (\ref{3.20e})-(\ref{3.25e}) and
therefore of Proposition 3.2. \par \nobreak \hfill $\sq$ \par

\noi {\bf Remark 3.2.} The parameters $\beta$, $\ell$ $\lambda_0$,
$k_+$ and $\mu$ play the same role and basically satisfy the same
conditions as in II. The parameter $\beta$ fixes the splitting of $B$
into long and short range parts and thereby fixes the auxiliary system
(1.20). The parameters $\ell$ and $\lambda_0$ fix the function space
where that system is solved. The regularity of $w$ is already fixed,
$\ell$ fixes the regularity of $s$ and $\lambda_0$ fixes the rates of
convergence in time. The parameters $k_+$ and $\mu$ fix the regularity
of $w_+$ and the vanishing of $(FA_+, F\dot{A}_+)$ at $\xi = 0$, which
have to be sufficient, as expressed by the lower bounds (3.67) (3.68).
The condition (3.68) is the combination of (3.84) (3.92), which
coincide with (II.3.81) and (II.3.82) respectively, with equality
allowed and with $\beta_0 = \beta$.\\

\noi {\bf Remark 3.3.} The conditions on the parameters become simpler
in the special case $\beta = 1/3$ which optimizes (3.68). If one takes
in addition $3/2 < \ell \leq 2$, the remaining conditions reduce to
$$1 < \lambda_0 < 3/2 \quad , \quad k_+ \geq \lambda_0 + 2 \quad , \quad \mu \geq (3/2) \lambda_0 - 7/4\ (> - 1/4) \ . \eqno(3.94)$$
\vskip 5 truemm

\noi {\bf Remark 3.4.} The condition $k_+ \geq \lambda_0 + 2$ is used
only to estimate the term $\nabla h \cdot \partial_t \nabla w_0$ in the
estimate of $\partial_t R_1$ in $L^2$ (see (3.87)). Everywhere else the
condition $k_+ \geq 3$ is sufficient. The latter condition would also
be sufficient for that term if one were using a better estimate of
$\parallel \nabla h \parallel_{\infty}$ than follows from (3.63) and
Sobolev inequalities (see Proposition 7.4 in \cite{4r} for estimates of
this type). \\

We finally comment briefly on the condition (3.61)$_\mu$ which
restricts the behaviour of $(FA_+, F\dot{A}_+)$ at $\xi = 0$. That
condition can be ensured by assuming sufficient decay of $(A_+,
\dot{A}_+)$ at infinity in space, possibly supplemented by some moment
conditions. We refer for details to the discussion at the end of
Section~II.3, from which we extract the following minor variation of
Lemma II.3.5, which is typical of the situation.\\

\noi {\bf Lemma 3.5.} {\it Let $-1/2 \leq  \mu < 1$. Let $(A_+, \dot{A}_+)$ satisfy (3.60) 
and in addition
$$x^2\ A_+ \in L^{2\vee 3/(\mu + 1)}\quad , \quad x^2 \dot{A}_+ \ , \ xA_+ \in L^{3/(2 + \mu )}$$
$$<x>^{1 + \mu + \varepsilon}  \dot{A}_+ \in L^1 \quad , \int \dot{A}_+ \ dx = 0\ ,$$
$$A_+, x \dot{A}_+ \in L^{3/(3 + \mu )} \qquad \hbox{\it for $\mu < 0$}\ ,$$
$$\int A_+\ dx = \int x \dot{A}_+ \ dx = 0 \quad , \quad <x>^{\mu + \varepsilon} A_+ \in L^1 \qquad
\hbox{for} \ \mu \geq 0 \ , $$

\noi for some $\varepsilon > 0$. Then $(3.61)_{\mu}$ holds.} 

\mysection{Wave operators and asymptotics for (u, A)} 
\hspace*{\parindent}
In this section we complete the construction of the wave operators for
the system (1.1) (1.2) and we derive asymptotic properties of solutions in their range. The
construction relies in an essential way on Proposition 3.2. So far we have worked with the
system (1.20) for $(w, s)$ and the first task is to reconstruct the phase $\varphi$.
Corresponding to $S = s_0 + s_1$, we define $\phi = \varphi_0 + \varphi_1$ where
\beq
\label{4.1e}
\varphi_0 =  \int^1_t dt'\ t'^{-1} \ B_L \left ( w_0(t'), w_0(t')\right ) \ ,
\eeq
\beq
\label{4.2e}
\varphi_1 = \int^t_0 dt' (1/2) |s_0(t')|^2 - 2 \int^t_0 dt' \ t'^{-1}
\ B_L \left ( w_0(t') , w_1(t')\right ) \ , 
\eeq

\noi so that $s_0 = \nabla \varphi_0$ and $s_1 = \nabla \varphi_1$. \par

Let now $(w, s)$ be the solution of the system (1.20) constructed in Proposition
3.1 and let $(q, \sigma ) =(w, s) - (W,S)$. We define\footnote{We take this opportunity to correct an omission in II, where the terms $B_{0L}$ and $B_L(W + W_1, w_2)$ are missing in (II.4.3). This has no incidence on the rest of II.}
$$\psi =  \int^t_0 dt' (1/2) \left ( \sigma \cdot (\sigma + 2S) + s_1 \cdot (
s_1 + 2s_0)\right ) (t')$$
 \beq
\label{4.3e}
- \int^t_0 dt' \ t'^{-1} \left ( B_{0L} + B_L (q, q) + 2
B_L(W,q) + B_L (w_1, w_1) + B_L (W + W_1, w_2) \right ) (t')
\eeq

\noi which is taylored to ensure that $\nabla \psi = \sigma$, given the fact that $s_0$, $s_1$
and $\sigma$ are gradients. The integral is easily seen to converge in $\dot{H}^1$ (see (I.8.4) and (3.91)-(3.93) with $m = 0$), and to
satisfy 
\beq
\label{4.4e}
\parallel \nabla \psi \parallel_2 \ = \ \parallel \sigma \parallel_2 \ \leq C\
t^{\lambda_0} \ . 
\eeq

\noi Finally we define $\varphi = \phi + \psi$ so that $\nabla \varphi = s$, and $(w,
\varphi )$ solves the system (1.18). For more details on the reconstruction of
$\varphi$ from $s$, we refer to Section 8 of I. \par

We can now define the wave operators for the system (1.1) (1.2) as follows. We start from the
asymptotic state $(u_+, A_+, \dot{A}_+)$ for $(u, A)$. We define $w_+ = Fu_+$, we define
$B_0$ by (1.4) (1.13), namely
$$A_0 = \dot{K}(t) \ A_+ + K(t) \ \dot{A}_+ = t^{-1} \ D_0(t)\ B_0(1/t) \ ,$$

\noi and we define $(W,S)$ by (3.42)-(3.47).\par

We next solve the system (1.20) with initial time zero by Proposition 3.2 and we reconstruct $\varphi$ from
$s$ as explained above, namely $\varphi = \varphi_0 + \varphi_1 + \psi$ with $\varphi_0$, $\varphi_1$ and $\psi$ defined by
(\ref{4.1e}) (\ref{4.2e}) (\ref{4.3e}) with $(q, \sigma ) = (w, s) - (W,S)$. We finally substitute $(w, \varphi )$ thereby
obtained into (1.11) (1.3) thereby obtaining a solution $(u, A)$ of the system (1.1) (1.2). The wave operator is defined as
the map $\Omega : (u_+, A_+, \dot{A}_+) \to (u, A)$. \par

We now turn to the study of the asymptotic properties of $(u, A)$ and in particular of its convergence to its asymptotic form $(u_a, A_a)$ defined in a natural way (compare with (1.3) (1.11)) by
\beq
\label{4.5e}
u_a(t) = M(t)\ D(t)\exp \left ( i\phi (1/t)\right ) \overline{W}(1/t)
\eeq
\beq
\label{4.6e}
A_a(t) = A_0(t) + A_1\left ( |D(t)\ W(1/t)|^2 \right ) .
\eeq

\noi In order to compare $u$ with $u_a$, we need some estimates of the
difference $\exp (- i \varphi ) w - \exp (-i \phi ) W$. The following
lemma is based on part of the estimates of Proposition 3.2 but does not
assume that $(w, \nabla \varphi )$ is a solution of the auxiliary
system (1.20). \\

\noi {\bf Lemma 4.1.} {\it Let $0 < \beta < 1$, $\lambda_0 > 1$, $0 <
\tau \leq 1$ and $I = (0, \tau ]$. Let $W \in {\cal C}(I, H^2) \cap
{\cal C}^1(I, H^1)$ satisfy (\ref{3.20e}) (\ref{3.21e}) and
\beq
\label{4.7e}
\parallel \Delta W \parallel_2 \ \leq a_1\ t^{-1/2}
\eeq

\noi for all $t \in I$. Let $(w, \nabla \psi ) \in {\cal C} (I, H^2
\oplus H^1) \cap {\cal C}^1(I, L^2 \oplus L^2)$ (with $\psi$,
$\partial_t \psi \in {\cal C}(I, L^6))$ satisfy the estimates
\beq
\label{4.8e}
\parallel w(t) - W(t) \parallel_2 \ \leq C\ t^{\lambda_0} 
\eeq
\beq
\label{4.9e}
\parallel \partial_t (w(t) - W(t)) \parallel_2 \ \vee \ \parallel \Delta (w(t) - W(t)) \parallel_2 \ \leq C\ t^{\lambda_0 - 1} 
\eeq
\beq
\label{4.10e}
\parallel \nabla^{m+1} \partial_t^j \psi (t) \parallel_2 \ \leq C\ t^{\lambda_0 - j - \beta m} \qquad \hbox{\it for $j, m = 0, 1\ , \ j+ m \leq 1$} 
\eeq

\noi for all $t \in I$.\par

(1) Let $f = \exp (- i \psi ) w - W$. Then the following estimates hold for all $t \in I$~:
\beq
\label{4.11e}
\parallel f(t) \parallel_2 \ \leq C\ t^{\lambda_0}\ , 
\eeq
\beq
\label{4.12e}
\parallel \partial_t f(t) \parallel_2 \ \vee \ \parallel \Delta f(t) \parallel_2 \ \leq C\ t^{\lambda_0 - 1}\ . 
\eeq

(2) Let $\phi \in {\cal C}(I, W_{\infty}^2) \cap {\cal C}^1 (I, L^{\infty})$ satisfy the estimates
$$\left \{ \begin{array}{l} \parallel \partial_t \phi\parallel_{\infty} \ \leq C\ t^{-1} \\ \\ \parallel \nabla \phi\parallel_{\infty}\ \vee \ \parallel \Delta \phi\parallel_{\infty}\ \leq C(1 - \ell n\ t) \end{array} \right . \eqno(4.13)$$

\noi for all $t \in I$. Let $\varphi = \phi + \psi$ and $g = \exp (- i \phi )f = \exp (-i \varphi ) w - \exp (- i \phi ) W$. Then the following estimates hold for all $t \in I$~:}
$$\parallel g(t) \parallel_2 \ \leq C\ t^{\lambda_0}\ , \eqno(4.14)$$
$$\parallel \partial_t g(t) \parallel_2 \  \vee \ \parallel \Delta g(t) \parallel_2 \ \leq C\ t^{\lambda_0 - 1}\ . \eqno(4.15)$$
\vskip 5 truemm

\noi {\bf Proof.} \underbar{Part (1)}. We write
$$f = \left ( \exp (- i \psi ) - 1 \right ) w + w - W$$

\noi and we estimate
$$\parallel f \parallel_2\ \leq \ \parallel \psi \parallel_6\ \parallel w \parallel_3 \ + \ \parallel w - W \parallel_2$$
$$\leq \ C\parallel \nabla \psi \parallel_2\ \parallel w \parallel_3 \ + \ \parallel  w - W\parallel_2\ \leq C\ t^{\lambda_0} \ .$$

\noi Next
$$\partial_t f = - i \exp (-i \psi ) \left ( \partial_t \psi \right ) w + \exp (- i \psi ) \partial_t (w-W) + \left ( \exp (-i \psi ) - 1 \right ) \partial_t W$$

\noi so that
$$\parallel \partial_t f \parallel_2\ \leq \ \parallel \partial_t \psi \parallel_6\ \parallel w \parallel_3 \ + \ \parallel \partial_t (w - W) \parallel_2 \ + \ \parallel \psi \parallel_6\ \parallel \partial_t W \parallel_3 \ $$
$$\leq C\ \left ( t^{\lambda_0 - 1} + t^{\lambda_0 - 1/2}\right ) \leq C\ t^{\lambda_0 - 1}$$

\noi by (\ref{3.20e}) (\ref{3.21e}) and (\ref{4.8e})-(\ref{4.10e}). \par

Finally
$$\Delta f = - \exp (-i \psi ) \left ( i \Delta \psi + |\nabla \psi |^2\right ) w - 2 i \exp (-i \psi ) \nabla \psi\cdot \nabla w$$
$$+ \exp (- i \psi ) \Delta (w-W) + \left ( \exp (-i \psi ) - 1 \right ) \Delta W$$

\noi so that
$$\parallel \Delta f \parallel_2\ \leq \ \left ( \parallel \Delta  \psi \parallel_2\ + \ \parallel \nabla \psi \parallel_4^2 \right )  \ \parallel w \parallel_{\infty} \ + \ 2\parallel \nabla \psi \parallel_6 \ \parallel \nabla w \parallel_3$$
$$+\ \parallel \Delta (w - W) \parallel_2 \ +\ \parallel \psi \parallel_{\infty} \ \parallel \Delta W \parallel_2$$
$$\leq C\ \left ( t^{\lambda_0 - \beta} + t^{2\lambda_0 - 3 \beta /2} + t^{\lambda_0 -1} + t^{\lambda_0-(\beta + 1)/2}\right ) \leq C\ t^{\lambda_0 - 1}$$

\noi by (\ref{3.20e}) and (\ref{4.7e})-(\ref{4.10e}).\\

\noi \underbar{Part 2}. (4.14) is obvious. Next
$$\parallel \partial_t g \parallel_2\ \leq \ \parallel \partial_t \phi \parallel_{\infty}\ \parallel f \parallel_2 \ + \ \parallel \partial_t f \parallel_2 \ \leq C\ t^{\lambda_0 - 1}$$
$$\parallel \Delta g \parallel_2\ \leq \ \left ( \parallel \Delta  \phi \parallel_{\infty}\ + \ \parallel \nabla \phi \parallel_{\infty}^2 \right )  \ \parallel f \parallel_2\ + \ 2\parallel \nabla \phi \parallel_{\infty} \ \parallel \nabla f \parallel_2\ + \ \parallel \Delta f \parallel_2$$
$$\leq C\ \left ( t^{\lambda_0} (1 - \ell n\ t)^2 + t^{\lambda_0 - 1/2} (1 - \ell n\ t) + t^{\lambda_0-1} \right ) \leq C\ t^{\lambda_0 - 1}\ . $$
\nobreak \hfill $\sq$ \par

In order to state the asymptotic properties of $u$, it is convenient to use the related function (compare with (1.7) (1.11))
$$\widetilde{u}(t) = U(-t) \ u(t) = M(t)^* F^* \exp \left ( i \varphi (1/t)\right ) \overline{w}(1/t) \eqno(4.16)$$

\noi and its asymptotic form
$$\widetilde{u}_a(t) = U(-t) \ u_a(t) = M(t)^* F^* \exp \left ( i \phi (1/t)\right ) \overline{W}(1/t) \ . \eqno(4.17)$$

We can now state the asymptotic properties of $(u, A)$.\\

\noi {\bf Proposition 4.1.} {\it Let $0 < \beta < 3/5$, $\ell > 3/2$
and $1 \vee \beta (\ell + 1) < \lambda_0 < 3/2$. Let $k_+$ and $\mu$
satisfy (3.67) (3.68), let $w_+ = Fu_+ \in H^{k_+}$, let $(A_+,
\dot{A}_+)$ satisfy (2.16) (3.60) (3.61)$_{\mu}$. Let $(W,S)$ be
defined by (3.42) (3.47). Let $(w,s)$ be the solution of the system
(1.20) in $(0, \tau ]$ obtained in Proposition 3.2 and let $I = [T,
\infty )$ where $T = \tau^{-1}$. Let $\phi = \varphi_0 + \varphi_1$ and
$\varphi = \phi + \psi$ be defined by (\ref{4.1e})-(\ref{4.3e}). Define
$(u, A)$ by (1.3)-(1.5) (1.11), define $(u_a, A_a)$ by
(\ref{4.5e}) (\ref{4.6e}) and $(\widetilde{u}, \widetilde{u}_a)$ by
(4.16) (4.17). Then $\widetilde{u} \in {\cal C}(I, FH^2) \cap {\cal
C}^1(I, L^2)$, $A \in {\cal C}(I, H^1) \cap {\cal C}^1(I, L^2)$, $(u,
A)$ solves the system (1.1) (1.2) in $I$ and $(u, A)$ behaves
asymptotically in time as $(u_a, A_a)$ in the sense that the following
estimates hold for all $t \in I$~:
$$\parallel \partial_t^j  (\widetilde{u} (t) - \widetilde{u}_a(t)) \parallel_2\ \leq \ C\left ( a_+, A_+, \dot{A}_+\right ) t^{-\lambda_0 - j}\qquad {\it for}\ j = 0,1\ , \eqno(4.18)$$
$$\parallel x^2  (\widetilde{u} (t) - \widetilde{u}_a(t)) \parallel_2\ \leq \ C\left ( a_+, A_+, \dot{A}_+\right ) t^{-\lambda_0 +1}\ ,\eqno(4.19)$$
$$\parallel u (t) - u_a(t) \parallel_r\ \leq \ C\left ( a_+, A_+, \dot{A}_+\right ) t^{-\lambda_0 + \delta (r)/2}\qquad {\it for}\ 2 \leq r \leq \infty   \ . \eqno(4.20)$$

\noi Furthermore $A - A_a \in {\cal C}(I, H^3)$ and the following
estimates hold for all $t\in I$~:
$$\parallel A (t) - A_a(t) \parallel_2\ \leq \ C\left ( a_+, A_+, \dot{A}_+\right ) t^{-\lambda_0 + 1/2}\ , \eqno(4.21)$$
$$\parallel \omega^{m+1} \left ( A (t) - A_a(t)\right )  \parallel_2\ \leq \ C\left ( a_+, A_+, \dot{A}_+\right ) t^{-\lambda_0 - 1/2-m/2}\left ( 1 + t^{m/2 - 3/4}\right )  \eqno(4.22)$$

\noi for $0 \leq m \leq 2$. Here $a_+ = |w_+|_{k_+}$ and the constants
$C(a_+, A_+, \dot{A}_+)$ depend on $(A_+, \dot{A}_+)$ through the norms
associated with (2.16) (3.60) (3.61)$_{\mu}$.}\\

\noi {\bf Proof.} The regularity of $u$ follows immediately from that
of $(w, \varphi )$ implied by Proposition 3.2, from the definition of
$\varphi$, from (4.16) and from the commutation relations
 $$\left \{ \begin{array}{l} x^2 M^*F^* = - M^* F^* \Delta\\ \\ i \partial_t M^* F^* = M^* F^* \left ( i \partial_t + \left ( 2t^2\right )^{-1} \Delta \right ) \ .\end{array} \right . \eqno(4.23)$$
 
\noi The estimates (4.18) (4.19) follow from Lemma 4.1, more precisely
from (4.14) (4.15) and from (4.23) again, once we have verified the
assumptions of that lemma. Now (\ref{4.7e}) follows from Lemma 3.3,
especially (3.52) (3.55) and from (3.74) with $m = 2$. The estimates
(\ref{4.8e}) (\ref{4.9e}) are a rewriting of (3.69) (3.70), while
(\ref{4.10e}) is a special case of (3.71) (3.72) with $m = 0, 1$.
Finally (4.13) follows  from the definition of $\phi$ and from Lemma
3.3, which implies that 
 $$\parallel \omega^{m+1} \ \phi \parallel_2\ \leq \ C(1- \ell n\ t) \eqno(4.24)$$
 $$\parallel \omega^{m+1} \ \partial_t \phi \parallel_2\ \leq \ Ct^{-1}$$
 
 \noi for $0 \leq m \leq 3 \leq k_+$.\par
 
 The estimate (4.20) follows from the fact that
 $$ u(t) - u_a(t) = M(t)\ D(t)\ \overline{g}(1/t) \eqno(4.25)$$
 
 \noi where $g$ is defined in Lemma 4.1, so that for $2 \leq r \leq \infty$
 $$\parallel u(t) - u_a(t) \parallel_r\ = \ t^{-\delta (r)} \parallel g(1/t)\parallel_r$$
 $$\leq \ C\ t^{-\delta (r)} \parallel g(1/t)\parallel_2^{1 - \delta (r)/2}\ \parallel \Delta g(1/t) \parallel_2^{\delta (r)/2} \ \leq C\ t^{-\lambda_0 - \delta (r) /2} \eqno(4.26)$$
 
 \noi by (4.14) (4.15).\par
 
 The regularity properties of $A$ follow immediately from the
assumptions on $(A_+, \dot{A}_+)$, from (1.13)-(1.15) and from the
regularity of $w$ implied by Proposition 3.2. In order to derive the
estimates (4.21) (4.22), we note that
$$A(t) - A_a(t) = t^{-1}\ D_0(t)\ B_1(q, q+2W)(1/t) \eqno(4.27)$$
 
 \noi with $q = w - W$. We estimate 
 $$\parallel B_1(q, q+2W)(t) \parallel_2\ \leq \ C\ I_{-1} \left ( \parallel q\parallel_2\ \parallel q + 2W \parallel_3\right ) \ \leq\ C\ t^{\lambda_0} \eqno(4.28)$$
 
 \noi by (2.8) and (3.69),
 $$\parallel \omega^{m+1}\ B_1(q, q+2W)(t) \parallel_2\ \leq\ C\ I_m \Big ( \parallel \omega^m q \parallel_2\ \parallel q + 2W \parallel_{\infty}$$
 $$+\ \parallel \omega^{\delta} q \parallel_2\ \parallel \omega^m\ w_0 \parallel_{3/\delta}\ + \ \parallel q \parallel_{\infty}\left ( \parallel \omega^m w_1 \parallel_2\ + \ \parallel \omega^m \ w_2 \parallel_2\right ) \Big )$$
 
\noi by (2.8) and Lemma 2.2, with $0 < \delta < 3/2$, 
$$\cdots \leq C \left ( t^{\lambda_0 - m/2} + t^{\lambda_0 - \delta /2} + t^{\lambda_0 - 3/4} \left ( t (1 - \ell n\ t) + t^{3/2-m}\ 1_m \right ) \right )$$

\noi for $0 \leq m \leq 2$ and $\delta \geq m - 3/2$, by (3.69) (3.70) (3.52) (3.55) (3.74),
$$\cdots \leq C\ t^{\lambda_0 - m/2} \left ( 1 + t^{3/4 - m/2}\right )\ . \eqno(4.29)$$

\noi Combining (4.28) (4.29) with $t$ changed into $1/t$ and (4.27)
yields (4.21) (4.22).\par \nobreak \hfill $\sq$\par

We conclude this section with some comments on the relation between
this paper and I, II. In all cases the regularity of $w$ includes $w\in
{\cal C}(H^k)$, with $1 < k < 2$ and a support condition on $Fu_+$ in
$I$, and with $1 < k < 3/2$ and no support condition in II. The present
paper covers the case $k=2$ with no support condition. In II, the
correction term $w_2$ allowing to eliminate the support condition could
be omitted in the statement of the final result (see Proposition
II.4.1, part (3)) because $w_2$ satisfied the same estimates as $q$
(namely (II.4.8) (II.4.9) with $q = w- W$). This is no longer the case
here, and $w_2$ has to be kept in the final result. The reason is that
whereas $\parallel w_2 \parallel_2$ and $\parallel \partial_t
w_2\parallel_2$ are estimated as $\parallel q\parallel_2$ and
$\parallel \partial_t q\parallel_2$, we have only $\parallel \Delta
w_2\parallel_2\ \leq C t^{-1/2}$ by (3.74), as compared with $\parallel
\Delta q\parallel_2 \ \leq C t^{\lambda_0 - 1}$.\par

We finally mention two minor differences in the statement of the
asymptotic properties of $u$. In I, II we have used the operator $J = x
+ it \nabla$ and norms of the type $\parallel |J(t)|^k u\parallel_2$.
By the commutation relation

$$J(t) = U(t) \ x \ U(-t)$$

\noi this is equivalent to $\parallel |x|^k \widetilde{u}\parallel_2$.
We have used $\widetilde{u}$ in the present paper because it allows for
simpler statements when considering time derivatives. The second minor
difference is that in I, II we have stated the convergence of $u$ to
$u_a$ in terms of $\exp (i \phi ) (u-u_a)$ instead of $u-u_a$. We could
as well have stated them also in terms of $u - u_a$ by using a suitable
variant of Lemma 4.1, part (3). This possibility arises mainly because
the WS system belongs to the borderline long range case, where the
correcting phase $\phi$ is only logarithmic.

\end{document}